\documentclass[11pt]{article}
\usepackage{graphicx,mathrsfs,amssymb}
\usepackage{amsmath,latexsym,amsbsy, color}
\usepackage[linkcolor=blue,colorlinks=true,citecolor=red,linktocpage]{hyperref}
\usepackage{graphicx}

\def\ds{\displaystyle}
\def\forall{\hbox{for all}~}

\def\ve{\varepsilon}
\def\n{\noindent}

\def\R{\mathbb{R}}

\def\vp{\varphi}

\def\vs{\vskip 2em}

\def\v{\vskip 1em}
\def\O{{\cal O}}

\def\X{{\cal X}}

\def\C{{\cal C}}
\def\bega{\begin{array}}
\def\enda{\end{array}}
\def\begi{\begin{itemize}}
\def\endi{\end{itemize}}
\def\ov{\overline}
\def\Tilde{\widetilde}

\def\bel{\begin{equation}\label}
\def\eeq{\end{equation}}
\def\sqr#1#2{\vbox{\hrule height .#2pt
\hbox{\vrule width .#2pt height #1pt \kern #1pt
\vrule width .#2pt}\hrule height .#2pt }}
\def\square{\sqr74}
\def\endproof{\hphantom{MM}\hfill\llap{$\square$}\goodbreak}

\newtheorem{theorem}{Theorem}[section]

\newtheorem{lemma}{Lemma}[section]

\newtheorem{remark}{Remark}[section]
\newtheorem{definition}{Definition}[section]

\usepackage{marginnote}

\begin{document}
\title{\bf Piecewise regular solutions to scalar balance laws with singular nonlocal sources}\vs
\author{\it Lorena Bociu$^{(1)}$, Evangelia Ftaka$^{(1)}$, Khai T. Nguyen$^{(1)}$,  and Jacopo Schino$^{(2)}$\\
\\
		{\small $^{(1)}$ Department of Mathematics, North Carolina State University}\\
		{\small $^{(2)}$ Faculty of Mathematics, Informatics and Mechanics, University of Warsaw}\\
		{\small e-mails: ~lvbociu@ncsu.edu, ~eftaka@ncsu.edu,   ~ khai@math.ncsu.edu, ~j.schino2@uw.edu.pl }}
\maketitle
\begin{abstract} The present paper  establishes a local well-posed result for  piecewise regular solutions with single shock of  scalar balance laws with  singular integral of convolution type kernels.  In a neighborhood of the shock curve, a detailed description of the
solution is provided for a general class of
initial data.
		\medskip
		
		\n {\bf AMS Mathematics Subject Classification.} 35B65, 76B15.
	\end{abstract}

%\begin{center}
%	\begin{minipage}{11.5cm}
%		\tableofcontents
%	\end{minipage}
%\end{center}

%\tableofcontents

\section{Introduction}
\label{sec:1}
\setcounter{equation}{0}
\setcounter{equation}{0}
We consider a scalar balance law in one space dimension with  a singular source term
\bel{BL}
u_t+f(u)_x~=~{\bf G}[u],
\eeq
where  $u:[0,\infty[\times\R\to\R$ is the state variable, $f:\R\to\R$ is a $C^4$ strictly convex flux, i.e.,
\bel{strict-conv}
\theta\cdot f(x_1)+(1-\theta)\cdot f(x_2)~>~f(\theta\cdot x_1+(1-\theta)\cdot x_2),\quad \theta\in ]0,1[, \ \  x_2\neq x_1,
\eeq
and  ${\bf G}$ is a singular integral of convolution type defined by convolution with a kernel $K$ that is locally integrable on $\R\backslash\{0\}$, in the sense that
\bel{G}
{\bf G}[g](x)~=~\lim_{\ve\to 0+}\int_{|y-x|>\ve}K(x-y)\cdot g(y)~dy.
\eeq
Here we work under the following assumptions on $K$, as it is typically done in applications:
\begin{itemize}
\item [\bf (H1)]  The kernel $K\in \mathcal{C}^2(\R\backslash \{0\})$ takes the form of 
$$K~=~K_1+K_2\qquad \mathrm{with}\qquad K_2\in{\bf L}^1(\R),$$
and the singular part $K_1$  is odd;
%\[
%K=K_1+K_2,\qquad \text{with} \  K_2~\in~{\bf L}^1(\R),
%\]
%and the singular part $K_1$ is odd.
\item [{\bf (H2)}] There exist a constant $C>0$ such that 
\bel{K1}
\left|K^{(i)}(x)\right|~\leq~{C\over |x|^{i+1}}\qquad\forall i=0,1,2.
\eeq
\end{itemize}
The conditions {\bf (H1)}-{\bf (H2)} ensure that  the Fourier transform of $K$  is essentially bounded. Thus,   ${\bf G}:{\bf L}^{2}(\R)\to {\bf L}^2(\R)$ is a bounded operator \cite{S}, i.e., 
\[
\|{\bf G}(g)\|_{{\bf L}^{2}(\R)}~\leq~\|{\bf G}\|_{\infty}\cdot \|g\|_{{\bf L}^2(\R)},\qquad g\in {\bf L}^{2}(\R).
\]
%
%Here, as it is typically done in applications, the kernel $K$ takes the form of 
%\[
%K=K_1+K_2,\qquad \text{with} \  K_2~\in~{\bf L}^1(\R),
%\]
%
%
%We work under the following assumptions on $\bf G$, as it is typically done in applications: 
%\begin{itemize}
%\item [\bf (H1)] ${\bf G}:{\bf L}^{2}(\R)\to {\bf L}^2(\R)$ is a bounded operator;
%\item [{\bf (H2)}] The kernel $K$ is in $\mathcal{C}^2(\R\backslash \{0\})$ and satisfies
%\bel{K1}
%\left|K^{(i)}(x)\right|~\leq~{C\over |x|^{i+1}}\qquad\forall i=0,1,2.
%\eeq
%\end{itemize}
%Under the assumption {\bf (H2)}, the ${\bf L}^{\infty}$ bound on the Fourier transform of $K$ ensures the continuity condition {\bf (H1)} (see e.g in \cite{S}). In particular, assumption {\bf (H1)} holds if  the kernel $K$ takes the form of 
%\[
%K=K_1+K_2,\qquad \text{with} \  K_2~\in~{\bf L}^1(\R),
%\]
%and the singular part $K_1$ is odd and satisfies {\bf(H2)}. 

Equation (\ref{BL}) has  an interesting structure since the scalar conservation law generates a contractive semigroup on ${\bf L}^1(\R)$, but the operator ${\bf G}$ may be  discontinuous and unbounded as an operator on ${\bf L}^1(\R)$. In the archetypal case when  
$$
f(u)~=~{u^2\over 2},\qquad\quad K(x)~=~\ds {1\over \pi x},
$$
equation (\ref{BL}) is well-known as the Burger-Hilbert equation which was introduced by Biello and Hunter in \cite{BiH} as a model for surface waves with constant frequency.  A lower bound on the maximal time of existence for smooth solutions was studied in \cite{BiH, HI, HIDT}, the formation of singularities  \cite{CST} and the local asymptotic behavior of a solution up to the time when a new shock is formed in finite time was investigated in \cite{Y}, and the global existence of entropy weak  solutions  was proved in \cite{BN}, together with a partial uniqueness result.  Recently, piecewise regular solutions with a single shock for the ``well-prepared" initial data have been constructed  in \cite{BZ,KV}. To complete an asymptotic description of a solution to the Burgers-Hilbert equation in a neighborhood of a point $y_0$ where two shocks interact in \cite{BNK},
the result was extended to a bigger class of initial data 
\[
u(0,x)~=~\overline{w}(x-y_0)+\Big(c_1\cdot \chi_{\strut ]-\infty,y_0[}+c_2\cdot \chi_{\strut ]y_0,+\infty[}\Big)\cdot\psi (x-y_0),
\]
for some $\overline{w}\in\ H^2(\R\backslash\{0\})$, constants $c_1,c_2\in\R$, and $\psi(x)\in \mathcal{C}^{\infty}(\R\backslash\{0\})$ being a fixed even function with compact support,  smooth outside the origin and satisfying
\[\psi(x)~=~{2\over \pi}\cdot |x|\ln |x|\qquad\forall~~|x|\leq1.
\]
In the present article, we study  the unique piecewise regular solution with a single shock of (\ref{BL}) with a general class of initial data of the form
\bel{IC} u(0,x)~=~ \ov w(x-y_0)+\bar v(x-y_0),\qquad \ov w~\in ~H^2(\R\backslash\{0\}),
\eeq
with for some $3/4<\alpha< 1$,  
\bel{Xadef}\bega{l} 
\bar v(x)~\in~\X_\alpha ~\doteq~
\Bigg\{ v\in \C_c^0(\R)\cap\mathcal{C}^4(\R\backslash\{0\}):%~\ds ~\sup_{x\not= 0}    {|x|^{-\alpha} \over 1+|x|^{-\alpha}} 
%%\bigl|v(x)\bigr|~<~+\infty, \\[4mm]
%\qquad\qquad\qquad\qquad\qquad \ds\qquad 
\ds \sup_{x\not= 0}    {|x|^{i-\alpha} \over 1+|x|^{i-\alpha}}
\bigl|v^{(i)}(x)\bigr|<+\infty,\quad 1\leq i\leq 4\Bigg\}.\enda\eeq
Intuitively, $\X_\alpha$ is the space of functions that have an arbitrarily large derivative as $x\to 0\pm$ and grow like $|x|^{\alpha}$. However, this does not lead to the formation of additional shocks. Indeed, one expects that characteristics fall into the large shock before the blow-up of the gradient.  As in \cite{BZ,BNK}, the solutions are more regular than the usual weak entropy solutions and can be determined by integrating along characteristics. These correspond to the ``broad solutions" considered in \cite{Bbook, RY}. 

Below we recall the notion of piecewise regular solutions for scalar balance laws. 

\begin{definition}\label{d:12}
 {\it A function $u:[0,T]\times \R\to\R$ is called a {\bf piecewise regular solution} of (\ref{BL})
% 
% n entropy weak solution $u=u(t,x)$   of (\ref{BL}), defined 
%on the interval $t\in [0,T]$, will be called a {\bf piecewise regular solution}
 if there exist
 finitely many shock curves $y_1(t), \ldots, y_n(t)$ such that the following holds.
\begi
\item[(i)] The map $t\mapsto u(t,\cdot)\in H^1( \R\backslash\{y_1(t),\ldots, y_n(t)\})\cap H^2_{loc}( \R\backslash\{y_1(t),\ldots, y_n(t)\})$ satisfies 
\[
\sup_{t\in [0,T]}\left(\|u(t,\cdot)\|_{H^1( \R\backslash\{y_1(t),\ldots, y_n(t)\})}+\|u(t,\cdot)\|_{H^2( \R\backslash (\bigcup_{i=1}^{n}[y_i(t)-\delta,y_i(t)+\delta]))}\right)<\infty
\]
for every $\delta>0$ sufficiently small.
%
%
%For every $t\in [0,T]$, the norm  $\|u(t,\cdot)\|_{H^1( \R\backslash\{y_1(t),\ldots, y_n(t)\})}$ remains uniformly bounded. Moreover, for every $\delta>0$ sufficiently small, it holds  
%\[
%|u_{x}(t,x)|~\leq~M_{\delta}\qquad\forall (t,x)\in [0,T]\times\left[ \R\backslash \left(\bigcup_{i=1}^{n}(y_i(t)-\delta, y_i(t)+\delta)\right)\right].
%\]
\item[(ii)]  For each $i=1,\ldots,n$, the Rankine-Hugoniot conditions hold:
\bel{RH1} u_i^-(t)~\doteq~u(t, y_i(t)-)~>~u(t, y_i(t)+) ~\doteq~u_i^+(t),\eeq
\bel{RH2} \dot y_i(t)~=~ {f(u_i^-(t)) -f(u_i^+(t))\over u_i^-(t)-u_i^+(t)}\,.\eeq
\item[(iii)]  Along every characteristic curve $t\mapsto x(t)$ such that $\dot x(t)=f'(u(t,x(t)))$, 
one has
\bel{cc2} {d\over dt} u\bigl(t, x(t)\bigr) ~=~{\bf G}[u](x(t)).\eeq
\endi
}
\end{definition}
We note that in the above definition, as well as throughout the sequel, the upper dot denotes a derivative with respect to time.

The remainder of the paper is structured as follows.  Our main theorem is presented in Section \ref{sec:main}, along with the main steps of the proof. Section 3 develops various key a priori estimates on the source term which are necessary in the remaining steps of the proof. 
In Section 4 we construct the unique local piecewise regular solution to [(\ref{BL}),(\ref{IC})] as the limit of a convergent sequence of approximations. We conclude with an Appendix which contains some basic estimates on the singular kernel, the corrector term,  and related functions. 
\section{Main result}
\label{sec:main}
\setcounter{equation}{0}
 We establish the local existence and uniqueness of a piecewise regular solution with a single shock  to  the general scalar balance law with nonlocal singular sources (\ref{BL}) for  a large class of initial data defined in (\ref{IC})-(\ref{Xadef}):
\begin{equation}\nonumber
\begin{cases}
u_t+f(u)_x~=~{\bf G}[u],\\[3mm]
u(0,x)~=~ \ov w(x-y_0)+\bar v(x-y_0),\ \  \text{with} \  \ov w~\in ~H^2(\R\backslash\{0\}), \ \text{and} \ \bar v~\in~\X_\alpha.
\end{cases}
\end{equation}

Our main theorem is presented below.
%the existence and uniqueness of a piecewise regular solution with single shock to the Cauchy problem (\ref{BL}) with a class of  initial data in (\ref{IC}), locally in time.
\begin{theorem}\label{t:main} Given $y_0\in\R$ and $\bar{v}\in \X_\alpha$ with $\ds \alpha\in (3/4,1)$, for every $\overline{w}\in H^2\bigl(\R\backslash \{0\}\bigr)$ such that $\overline{w}(0-)>\overline{w}(0+)$,  the Cauchy problem (\ref{BL}) with initial data (\ref{IC})  admits a unique piecewise regular solution $u$ with a single shock $y_1(\cdot)$ defined for $t\in [0,T]$, for some $T>0$ sufficiently small. Moreover,  the map $t\mapsto u(t,y_1(t)\pm)$ is locally Lipchitz and  satisfies 
\bel{Lp-u}
\left|\dot{u}(t,y_1(t)\pm)\right|~\leq~ 2\Gamma_1 t^{\alpha-1}~~~ a.e.~t\in [0,T],
\eeq
for some constant $\Gamma_1>0$.
\end{theorem}
\begin{remark}
The local existence and uniqueness result can be  extended to the case of solutions with finitely many non-interacting shocks. Moreover, our result can be applied to both Fornberg-Whitham equation \cite{FW} and Burgers-Possion equation \cite{FS,GN,KN}. 
\end{remark}
\begin{remark}
 We point out that
 the lower bound $3/4$ on the constant $\alpha$ in the definition of $\X_\alpha$ is somehow sharp within our analysis. The best lower bound on $\alpha$ remains an open question.
\end{remark}
The main steps in the proof of our main theorem are introduced below, while the details are provided in the subsequent sections. 
\subsection{New coordinate system for solutions with one shock}
The first step in the proof of our main theorem consists in transferring the equation (\ref{BL}) to a new coordinate system so that the location  of the shock of the constructed piecewise regular solution  always remains at the origin. The details for the change of coordinates are included below. 

Assume that $u$ is  a piecewise regular solution of the balance law (\ref{BL}), with one single shock. 
By the Rankine-Hugoniot condition in (\ref{RH2}), the location $y(t)$ of the shock at time $t$ satisfies  
\[
\dot{y}(t)~=~{f(u^{-}(t))-f(u^{+}(t))\over u^{-}(t)-u^{+}(t) },\qquad u^{\pm}(t)~=~\lim_{x\to y(t)\pm}u(t,x).
\]
 As in \cite{BZ,BNK},  we shift the space coordinate, by replacing $x$ with $x-y(t)$, so that in the new coordinate system
the shock is always located at the origin. In these new coordinates, the Cauchy problem [(\ref{BL}), (\ref{IC})]  becomes
\bel{BL-1}
u_t+\left(f'(u)- {f(u^-(t))-f(u^+(t))\over u^{-}(t)-u^{+}(t)}\right)\cdot u_x~=~{\bf G}[u],
\eeq
with  initial data of the form 
\bel{IC1} u(0,x)~=~ \ov w(x)+\bar v(x),\qquad \ov w\in H^2(\R\backslash\{0\}),\qquad \bar{v}\in \X_\alpha.
\eeq
%\textcolor{blue}
%{
Let $\eta\in C^{\infty}(\R)$ be  an even cut-off function  which is nonincreasing on $[0,\infty[$ and satisfies
 \bel{eta}
\mathrm{supp}(\eta)~\subseteq~[-2,2],\qquad \eta(x)~=~1\qquad\forall x\in [-1,1].
\eeq
Without loss of generality, we can assume $\bar{v}(0) = 0$ with $\mathrm{supp}(v)\subseteq [-2,2]$. This assumption is justified due to the fact that  \eqref{IC1} can be rearranged as
\begin{equation*}
\overline{w} + \bar{v} = \Bigl[\overline{w} + \bar{v}-\Bigl(\bar{v} - \bar{v}(0)\Bigr)\cdot\eta\Bigr] + \Bigl(\bar{v} - \bar{v}(0)\Bigr)\cdot\eta~\in~ H^2(\R\backslash\{0\})+\X_{\alpha}
\end{equation*}

%with $\eta\in C^{\infty}(\R)$ being  an even cut-off function  which is nonincreasing on $[0,\infty[$ and satisfies
% \bel{eta}
%\mathrm{supp}(\eta)~\subseteq~[-2,2],\qquad \eta(x)~=~1\qquad\forall x\in [-1,1].
%\eeq

\subsection{Structure of piecewise regular solution}
\label{subsec:3}

The key idea in proving a local existence result for the Cauchy problem  (\ref{BL-1})-(\ref{IC1}) is  to look for solutions of the form
\bel{usw} u(t,x)~=~w(t,x)+\vp^{(w)}(t,x),\eeq
where $w(t,\cdot)$ belongs to $ H^2(\R\backslash \{0\})$ for $t>0$ and the corrector term 
$\vp^{(w)}$  depends explicitly on time $t$ such that $\vp^{(w)}(t,0)=0$ for all $t\geq 0$,  the strength of the jump 
\bel{JS}\sigma^{(w)}(t)~
\doteq ~ w^-(t) -w^+(t),\qquad \text{where}\ \ w^{\pm}(t)\,\doteq\,w(t,0\pm),\eeq
and also on the difference between the speed of characteristic and the  speed of the shock curve
\bel{speed1}
b^{(w)}_\pm(t)~\doteq~b^{(\omega)}(t,0\pm),\qquad~~ b^{(w)}(t,x)~\doteq~f'(w(t,x))-{f(w^-(t))-f(w^+(t))\over w^{-}(t)-w^{+}(t)}\,.
\eeq
%
%
%\bel{speed1}
%\begin{cases}
%b^{(w)}_-(t)&\doteq~\ds f'(w^-(t))-{f(w^-(t))-f(w^+(t))\over w^{-}(t)-w^{+}(t)}\,,\\[4mm]
%b^{(w)}_+(t)&\doteq~\ds f'(w^+(t))-{f(w^-(t))-f(w^+(t))\over w^{-}(t)-w^{+}(t)}\,.
%\end{cases}
%\eeq
%In order to make an appropriate ansatz for the function $\vp^{(w)}$, we observe from (\ref{G}) and (\ref{K1}) that if $g \in H^1(\R\backslash\{0\})$ with $\mathrm{supp}(g)\subseteq [-2,2]$, then for every $x\ne0$ we have 
%\bel{G1}
%{\bf G}[g](x)~=~\int_{-2}^{2}g'(y)\cdot \Lambda(x-y)~dy+\big[g(0+)-g(0-)\big]\cdot \Lambda(x).
%\eeq
%The function $\Lambda\in C^3(\R\backslash\{0\})$ is the  antiderivative of $K$, which is given by
Calling $\chi_{E}$ the indicator function of $E$, we define the function $\Lambda\in C^3(\R\backslash\{0\})$ as the  antiderivative of $K$, that is,
\bel{Lambda_1}
\Lambda(x)~=~\int_{0}^{x}K_1(y)dy+\Lambda_2(x),
\eeq
with $\Lambda_2$ being the even function as the  antiderivative of $K_2$ such that 
\bel{Lambda_2}
\Lambda_2(x)~=~
-\left(\ds\int_{x}^{2}K(y)dy\right)\cdot \chi_{]0,2]}(x)+\left(\ds\int_{2}^{x}K(y)dy\right)\cdot \chi_{]2,\infty[}(x),\  x\in ]0,\infty[.
\eeq
%
%\bel{Lambda_1}
%\Lambda(x)=
%-\left(\ds\int_{x}^{2}K(y)dy\right)\cdot \chi_{]0,2]}(x)+\left(\ds\int_{2}^{x}K(y)dy\right)\cdot \chi_{]2,\infty[}(x),\  x\in ]0,\infty[, \ \ \text{and}
%\eeq
%\bel{Lambda_2}
%\Lambda(x)=
%-\left(\ds\int_{x}^{-2}K(y)dy\right)\cdot \chi_{]-\infty,-2[}(x)+\left(\ds\int_{-2}^{x}K(y)dy\right)\cdot \chi_{]-2,0[}(x),\ x\in ]-\infty,0[.
%\eeq
In order to make an appropriate ansatz for the function $\vp^{(w)}$, we make use of the following property.
\begin{lemma}\label{G1}
If $g \in H^1(\R\backslash\{0\})$ with $\mathrm{supp}(g)\subseteq [-2,2]$, then for every $x\ne0$ we have
\begin{equation*}
{\bf G}[g](x)~=~\int_{-2}^{2}g'(y)\cdot \Lambda(x-y)~dy+\big[g(0+)-g(0-)\big]\cdot \Lambda(x).
\end{equation*}
\end{lemma}
{\bf Proof.} We consider the case $0 < x < 2$, the others being similar. Integrating by parts, we obtain
\[
\begin{split}
\int_{-2}^2 g'(y) &\cdot \Lambda(x-y) dy~=~\lim_{\ve\to 0+} \left[\int_{-2}^{-\ve} + \int_{\ve}^{x-\ve} + \int_{x+\ve}^2\right] g'(y) \cdot \Lambda(x-y) dy\\
&=~\lim_{\ve\to 0+} \left(\left[\int_{-2}^{-\ve} + \int_\ve^{x-\ve} + \int_{x+\ve}^2\right] g(y) \cdot K(x-y) dy+g(-\ve) \cdot \Lambda(x+\ve)-g(\ve) \cdot \Lambda(x-\ve)\right)\\
&\qquad\qquad\qquad\qquad\qquad +\lim_{\ve\to 0+} \left[ g(x-\ve) \cdot \Lambda(\ve)- g(x+\ve) \cdot \Lambda(-\ve)\right]\\
&=~ {\bf G}[g](x) - \big[g(0+) - g(0-)\big] \cdot \Lambda(x)+\lim_{\ve\to 0+} \left[ g(x-\ve) \cdot \Lambda(\ve)- g(x+\ve) \cdot \Lambda(-\ve)\right].
\end{split}
\]
From the assumption {\bf (H1)}-{\bf (H2)} and (\ref{Lambda_1}), one has 
\[
\lim_{\ve\to 0+}\ve\Lambda(\ve)~=~0,\qquad \lim_{\ve\to 0+} \big[\Lambda(\ve)-\Lambda(-\ve)\big]~=~\lim_{\ve\to 0+} \int_{-\ve}^{\ve}K_1(x)dx~=~0,
\]
and this yields 
\[
\lim_{\ve\to 0+} \left[ g(x-\ve) \cdot \Lambda(\ve)- g(x+\ve) \cdot \Lambda(-\ve)\right]~=~0.
\]
The proof is complete. 
%
%
%We consider the case $0 < x < 2$, the others being similar. Integrating by parts, we obtain
%\begin{align*}
%\int_{-2}^2 g'(y) \cdot \Lambda(x-y) dy & ~= \lim_{\ve\to 0+} \left[\int_{-2}^{-\ve} + \int_{\ve}^{x-\ve} + \int_{x+\ve}^2\right] g'(y) \cdot \Lambda(x-y) dy\\
%& ~=~ \lim_{\ve\to 0+} \left[\int_{-2}^{-\ve} + \int_\ve^{x-\ve} + \int_{x+\ve}^2\right] g(y) \cdot K(x-y) dy\\
%& \qquad + g(-\ve) \cdot \Lambda(x+\ve) + g(x-\ve) \cdot \Lambda(\ve) - g(\ve) \cdot \Lambda(x-\ve) - g(x+\ve) \cdot \Lambda(-\ve)\\
%& ~=~ {\bf G}[g](x) - \big[g(0+) - g(0-)\big] \cdot \Lambda(x)\\
%& \qquad - 2 \lim_{\ve\to 0+} \frac{g(x+\ve) \cdot \Lambda(-\ve) \cdot \ve - g(x-\ve) \cdot \Lambda(\ve) \cdot \ve}{2 \ve}
%\end{align*}
\endproof
Observe that Lemma \ref{G1} was implicitly used in \cite[formula (2.11)]{BZ} in the particular case where ${\bf G}$ is the Hilbert transform.
\medskip

At this point we note that equation  (\ref{BL-1}) can be approximated by the simpler equation
\bel{BL-2}
u_t+\left(f'(u)- {f(u^-(t))-f(u^+(t))\over u^{-}(t)-u^{+}(t)}\right)\cdot u_x~=~\big[u^{+}(t)-u^{-}(t)\big]\cdot \Lambda(x).
\eeq
Indeed, we expect that the solutions of (\ref{BL-1}) and (\ref{BL-2}) with the same initial data
will have the same asymptotic structure near the origin. Their difference will lie in  $H^2\bigl(\R\backslash\{0\}\bigr)$. With this in mind, let $\Phi:[0,\infty[\to \R$ be the  antiderivative of $\Lambda$ such that 
\bel{Phi}
\Phi(0)~=~0,\qquad \Phi(x)~=~ \int_{0}^{x}\Lambda(y)dy\qquad\forall x\in \R\backslash\{0\}.
\eeq
%and let $\eta\in C^{\infty}(\R)$ be  an even cut-off function  which is nonincreasing on $[0,\infty[$ and satisfies
% \bel{eta}
%\mathrm{supp}(\eta)~\subseteq~[-2,2],\qquad \eta(x)~=~1\qquad\forall x\in [-1,1].
%\eeq
Consider the function 
\bel{phi-1}
\phi(x,b)~\doteq~\eta(x)\cdot[\Phi(b)-\Phi(x+b)]\qquad\forall x,b\in \R.
\eeq
We  make the ansatz 
\bel{vp1}
\vp^{(w)}(t,x)~=~ \begin{cases}
\ds {\sigma^{(w)}(t) \over b^{(w)}_{-}(t)}\cdot \left[\phi(x,0)-\phi\left(x,-t\cdot b^{(w)}_-(t)\right)\right]\\[3mm]
\qquad\qquad + \eta(x)\cdot\left[\overline{v}\left(x-t\cdot b^{(w)}_-(t)\right)-\bar{v}\left(-t\cdot b^{(w)}_-(t)\right)\right]\qquad\mathrm{if}\quad x<0,\\[4mm]
\ds {\sigma^{(w)}(t)\over b^{(w)}_{+}(t)}\cdot \left[\phi(x,0)-\phi\left(x,-t\cdot b^{(w)}_+(t)\right)\right]\\[3mm]
\qquad\qquad +\eta(x)\cdot\left[\overline{v}\left(x-t\cdot b^{(w)}_+(t)\right)-\overline{v}\left(-t\cdot b^{(w)}_+(t)\right)\right]\qquad\mathrm{if}\quad x>0.
\end{cases}
\eeq
From (\ref{BL-1}), (\ref{usw}) and (\ref{speed1}), the remaining component $w(t,\cdot)$ from (\ref{usw}) solves the equation 
\bel{BH-2}
w_t+a(t,x,w)\cdot w_x~=~F(t,x,w),
\eeq 
where $a$ and $F$ are respectively given by
\bel{a}
a(t,x,w)~=~b^{(w)}(t,x)+f'\left(w+\vp^{(w)}\right)-f'(w),
\eeq
and
\begin{eqnarray}\label{F}
\hspace{.3in} F(t,x,w)&=&{\bf G}\left[\vp^{(w)}\right] (t,x) -\left[f'\left(w+\vp^{(w)}\right)-f'(w)\right]\cdot\vp^{(w)}_{x}(t,x)   \\  
& &\qquad\qquad +\left({\bf G}\left[w\right](t,x)-\left[ \vp^{(w)}_{t}(t,x)+b^{(w)}(t,x)\cdot \vp^{(w)}_{x}(t,x)\right]\right). \nonumber
\end{eqnarray}
\subsection{Construction of solution}
\label{subsec:4}
In order to finish the proof of Theorem \ref{t:main}, we construct solutions to the Cauchy problem (\ref{BH-2})-(\ref{F}) with initial data satisfying
\begin{equation}\label{ICw}
w(0,\cdot)~=~\ov w(\cdot)\in H^2(\R\backslash\{0\}),\qquad \ov w(0-)-\ov w(0+)~>~0.
\end{equation}
%{\color{blue}
%Observe that, from (\ref{IC1}), (\ref{usw}), (\ref{vp1}), and (\ref{ICw}), we have
%\[
%\overline{w}(x) + \overline{v}(x) = u(0,x) = w(0,x) + \varphi^{(w)}(0,x) = \overline{w}(x) + \eta(x) \overline{v}(x) = \overline{w}(x) + \overline{v}(x) - (1-\eta(x)) \overline{v}(x)
%\]
%and so, up to replacing $\overline{w}$ with $\overline{w} - (1-\eta) \overline{v} \in H^2(\R \setminus \{0\})$ in (\ref{ICw}), the formulas mentioned above are coherent.
%}
%\\
Following the  analysis in \cite{BZ},  the solution will be obtained as the limit of a sequence of approximations. 
Namely, consider a sequence of linear approximations constructed as follows.
As a first step, we set
\[ w_1(t,x)~=~\ov w(x)
\qquad\hbox{for all}~~t\geq 0, ~~x\in\R\,.\]
By induction, let $w_n$ be given. We define  $w_{n+1}$ to
be the solution of the semilinear, non-homogeneous  Cauchy problem
\bel{wne}
w_t + a(t,x,w_n)\cdot w_x~=~F(t,x,w),\qquad w(0,\cdot)~=~\ov w.\eeq 
The induction argument requires two main steps:
\begi
\item [(i).] Existence and uniqueness of solutions to each semilinear 
problem (\ref{wne}) such that 
\[
\sup_{n\geq 1}\left\{\sup_{t\in [0,T]}~\bigl\|w_n(t,\cdot)\bigr\|_{H^2(\R\setminus\{0\})}\right\}~<~\infty.
\]
\item [(ii).] Convergence  in the weaker norm $H^1(\R\backslash\{0\})$, which will follow from the contractive property below
\bel{conv1}
\sup_{t\in [0,T]}\bigl\|w_{n+1}(t,\cdot)-w_n(t,\cdot)\bigr\|_{H^1(\R\setminus\{0\})}< {1\over 2} \sup_{t\in [0,T]}\bigl\|w_{n}(t,\cdot)-w_{n-1}(t,\cdot)\bigr\|_{H^1(\R\setminus\{0\})}.
\eeq
\endi
The details for these three main steps will be provided in the following sections.

\section{Key estimates on the source term $F$}
\setcounter{equation}{0}
In order to prove existence and uniqueness of solutions to each linear 
problem (\ref{wne}) above, we need a priori estimates on the source term $F$ defined in (\ref{F}) and recalled here:
\begin{eqnarray}
F(t,x,w)&=& {\bf G}\left[\vp^{(w)}\right] (t,x) -\left[f'\left(w+\vp^{(w)}\right)-f'(w)\right]\cdot\vp^{(w)}_{x}(t,x) \nonumber \\  
& &\qquad\qquad +\left({\bf G}\left[w\right](t,x)-\left[ \vp^{(w)}_{t}(t,x)+b^{(w)}(t,x)\cdot \vp^{(w)}_{x}(t,x)\right]\right).  \nonumber
\end{eqnarray}
%Towards to the ,  we shall provide  in Lemma \ref{F1} some . In order to  do so,
First, we  rewrite the source $F$ in the following way: 
\begin{equation}\label{equivF}
F(t,x,w)~=~A^{(w)}(t,x)+B^{(w)}(t,x)-C^{(w)}(t,x), \ \ \text{with}
\end{equation} 
\bel{A}
A^{(w)}(t,x)~\doteq~{\bf G}\left[\vp^{(w)}\right] (t,x) -\left[f'\left(w+\vp^{(w)}\right)-f'(w)\right]\cdot\vp^{(w)}_{x}(t,x),
\eeq
\bel{B}
B^{(w)}(t,x)~\doteq~{\bf G}\left[w\right](t,x)-\sigma^{(\omega)}(t)\cdot \Lambda(x) \cdot \eta(x)\, \ \text{and}
%\left(f'\left(w\right)- {f(w^-(t))-f(w^+(t))\over w^{-}(t)-w^{+}(t)}\right)\cdot \psi^{(w)}_x(t,x),
\eeq
\bel{C}
C^{(\omega)}(t,x)~\doteq~\vp^{(w)}_{t}(t,x)+b^{(w)}(t,x) \cdot \vp^{(w)}_{x}(t,x)+\sigma^{(\omega)}(t)\cdot\Lambda(x)\cdot \eta(x).
\eeq

\subsection{Estimating the corrector term $\vp^{(w)}(t,x)$}
Given two constants $M_0,\delta_0>0$, we shall assume  that  $w(t,\cdot)$ is in $H^2(\R\backslash\{0\})$ and satisfies 
\bel{w-1-b}
\|w(t,\cdot)\|_{H^2(\R\backslash\{0\})}~\leq~M_0,\qquad  \sigma^{(w)}(t)~>~\delta_0\qquad\forall t\in [0,T]\,.
\eeq
In addition, the map $t\mapsto w^{\pm}(t)\doteq w(t,0\pm)$ is locally Lipschitz and
\bel{Lip-w}
\big|\dot{w}^{\pm}(t)\big|~\leq~\ell(t)\qquad a.e.~t\in ]0,T]
\eeq
for some function $\ell : \; ]0,T] \to \; ]0,\infty[$ such that $\ds\lim_{t \to 0+} \ell(t) = \infty$.
From (\ref{w-1-b}) and (\ref{JS}), we have the following estimates for all $x\in \R\backslash \{0\}$ and $t\in [0,T]$ 
\bel{sigma}
|w(t,x)|, |w_x(t,x)|~\leq~ 2M_0,\qquad~~ \big|\sigma^{(w)}(t)\big|~\leq~|w^+(t)|+|w^{-}(t)|~\leq~4M_0\,,
\eeq
and 
\bel{b-speed}
\begin{cases}
\big|b^{(w)}(t,x)\big|,\big|b^{(w)}_x(t,x)\big|~\leq~4\cdot\|f\|_{\C^2([-2M_0,2M_0])}M_0~\doteq~b_1,\\[4mm]
 \big|b^{(w)}_{xx}(t,x)\big|~\leq~\|f\|_{\C^3([-2M_0,2M_0])} \left(|w_{xx}(t,x)|+4 M_0^2\right).
\end{cases}
\eeq
%Moreover,
%\bel{bb}
%\left|b(t,x) - b_+^{(w)}(t)\right| \le 2\|f\|_{C^2([-2M_0,2M_0])} \cdot M_0 \cdot x.
%\eeq
By the strict convexity of $f$ in (\ref{strict-conv}), one has that
\bel{ul-b}
-b_1~\leq~b_{+}^{(w)}(t)~\leq~-b_{0}~<~0~<~b_{0}~\leq~b_{-}^{(w)}(t)~\leq~b_1\,,
\eeq
with  $b_0$ defined by 
\begin{equation}\label{b0def}
b_0~\doteq~\delta_0\cdot \min_{a,b\in [-2M_0,2M_0], b-a\geq \delta_0}\int_{0}^1\int_{0}^{1}f''(b-rs (b-a))sdrds~>~0\,.
\end{equation}
Moreover, assumption (\ref{Lip-w}) implies that  
\bel{b-deri}
\big|\dot{b}^{(w)}_{\pm}(t)\big|~\leq~2\|f\|_{\C^2([-2M_0,2M_0])}\cdot\ell(t)\qquad\forall t\in [0,T]\,.
\eeq
%
%there exist constants $b_0,b_1>0$ depending  on $\delta_0,M_0$ and $f$ such that 
%\bel{ul-b}
%-b_1~\leq~b_{+}^{(w)}(t)~\leq~-b_{0}~<~0~<~b_{0}~\leq~b_{-}^{(w)}(t)~\leq~b_1\,.
%\eeq
Using (\ref{sigma})-(\ref{ul-b}), we  provide some  estimates on the corrector term $\vp^{(w)}(t,x)$  defined  (\ref{Phi})-(\ref{vp1}). In the following, as usual, by the Landau symbol $\O(1)$ we shall denote a uniformly bounded quantity which does not depend on $M_0,\delta_0$ and $f$.
\begin{lemma}\label{vp-est} 
For every $0<|x|<1/4$ and $0<t<\ds {1\over  4b_1}$ with $b_1$ defined in (\ref{b-speed}), we have that  
\bel{vpb1}
\left|\vp^{(w)}(t,x)\right|~\leq~C_0\cdot |x|\left(\big|\ln |x|\big|+t^{\alpha-1}\right),
\eeq
\bel{vpb1_1}
\ds\left|{d\over dx}\vp^{(w)}(t,x)\right|~\leq~\ds C_0\cdot\left(\big|\ln |x|\big|+{1\over (|x| +t)^{1-\alpha}}\right)\,,
\eeq
\bel{vpb2}
\ds\left|{d^i\over dx^i}\vp^{(w)}(t,x)\right|\leq~\ds C_0\cdot\left({1\over |x|^{i-1}}+{1\over (|x| +t)^{i-\alpha}}\right),\qquad i\in\{2,3\}\,
\eeq
Moreover,  for $\delta>0$ small, we obtain that 
\[
\left\|\vp^{(\omega)}(t,\cdot)\right\|_{H^2(\R\backslash [-\delta,\delta])}~\leq~ {C_0}\cdot \delta^{\alpha-3/2},
\]
with  a constant $C_0>0$ satisfying
\[
\left\|\vp^{(w)}(t,\cdot)\right\|_{L^{\infty}}~\leq~C_0~\doteq~ \O(1)\cdot\left(M_0b_0^{-1}+1\right)\cdot \left(b_0^{\alpha-3}+1\right)\,.
\]
\end{lemma}
{\bf Proof.} We rewrite $\vp^{(w)}$ equivalently as $\vp^{(w)}=\Tilde{\vp}^{(w)}+\ov{\vp}^{(w)}$, with 

\bel{vp-tb}
\begin{cases}
\Tilde{\vp}^{(w)}(t,x)~=~\ds {\sigma^{(w)}(t) \over b^{(w)}_{-}(t)}\cdot \left[\phi(x,0)-\phi\left(x,-t\cdot b^{(w)}_-(t)\right)\right]\cdot\chi_{]-\infty,0[}\\
\qquad\qquad\qquad\qquad\qquad\qquad  +\ds {\sigma^{(w)}(t)\over b^{(w)}_{+}(t)}\cdot \left[\phi(x,0)-\phi\left(x,-t\cdot b^{(w)}_+(t)\right)\right]\cdot\chi_{]0,\infty[},\\[6mm]
\ov{\vp}^{(w)}(t,x)~=~\eta(x)\cdot\left[\overline{v}\left(x-t\cdot b^{(w)}_-(t)\right)-\bar{v}\left(-t\cdot b^{(w)}_-(t)\right)\right]\cdot \chi_{]-\infty,0[}\\
\qquad\qquad\qquad\qquad\qquad\qquad   +\eta(x)\cdot\left[\overline{v}\left(x-t\cdot b^{(w)}_+(t)\right)-\overline{v}\left(-t\cdot b^{(w)}_+(t)\right)\right]\cdot \chi_{]0,\infty[}.
\end{cases}
\eeq
From the assumption {\bf (H2)},  for all $x\in [-1/4,1/4]\backslash\{0\}$, it holds that
\[
|\Phi(x)|~\leq~\O(1)\cdot\big|x\ln|x|\big|,\quad |\Lambda(x)|~\leq~\O(1)\cdot \big|\ln|x|\big|,\quad K^{(i)}(x)~\leq~{C\over |x|^{i+1}},\quad i=0,1,2.
\]
Recalling (\ref{sigma}) and (\ref{ul-b}), we estimate  for every $x\in (0,1/4]$ and $0<t<\ds {1\over  4b_1}$ that 
\[
\begin{cases}
\big|\Tilde{\vp}^{(w)}(t,x)\big|~\leq~\ds{4M_0\over b_0}\cdot\left|\phi(x,0)-\phi\left(x,-t\cdot b^{(w)}_+(t)\right)\right|~\leq~\O(1)\cdot {M_0\over b_0}\cdot \big|x\ln |x|\big|,\\[4mm]
\left|\ds{d\over dx}\Tilde{\vp}^{(w)}(t,x)\right|~\leq~\ds{4M_0\over b_0}\cdot\left|\ds\Lambda\left(x-tb^{(w)}_+(t)\right)-\Lambda\left(x\right)\right|~\leq~\ds\O(1)\cdot {M_0\over b_0}\cdot \big|\ln x\big|,\\[4mm]
\left|\ds{d^i\over dx^i}\Tilde{\vp}^{(w)}(t,x)\right|~\leq~\ds\O(1)\cdot {M_0\over b_0}\cdot {1\over x^{(i-1)}},\qquad i\in \{2,3\},\\
%\ds{4M_0\over b_0}\cdot\left|\ds {d^{(i-1)}\over d^{(i-1)}x}\left[\Lambda\left(x-tb^{(w)}_+(t)\right)-\Lambda\left(x\right)\right]\right|~\leq~\ds\O(1)\cdot {M_0\over b_0}\cdot {1\over x^{(i-1)}},\quad i\in \{2,3\},\\
\end{cases}
\]
and 
\[
%\begin{cases}
\big|\ov{\vp}^{(w)}(t,x)\big|~=~\left|\overline{v}\left(x-t\cdot b^{(w)}_{+}(t)\right)-\overline{v}\left(-t\cdot b^{(w)}_{+}(t)\right)\right|~\leq~\ds\O(1)\cdot \min\left\{1,{x\over (b_0t)^{1-\alpha}}\right\},
\]
\[
\left|\ds{d^i\over dx^i}\ov{\vp}^{(w)}(t,x)\right|~=~\left|\ds{d^i\over dx^i}\overline{v}\left(x-t\cdot b^{(w)}_{+}(t)\right)\right|~\leq~\ds\O(1)\cdot {1\over (x + b_0t)^{i-\alpha}}\qquad i\in\{1,2,3\}.
%\end{cases}
\]
The same estimates  hold for  $x\in [-1/4,0)$, $0<t< 1/(4b_1)$, and this yields (\ref{vpb1})-(\ref{vpb2}).
\endproof

\subsection{Estimating source $F$}
The next  lemma  provides some  estimates on the function $F$ in (\ref{F}).  These estimates will be used in Lemma \ref{Pr0} to  establish a priori bounds on the approximate solutions of the linear, non-homogeneous  Cauchy problem (\ref{wne}).
\begin{lemma}\label{F1} Assume that $w:[0,T]\times\R\to\R$ satisfies assumptions (\ref{w-1-b})-(\ref{Lip-w}). Then  for $|x|<1/4$ and a.e. $t\in [0,T]$, we have
\bel{eqF1}
\begin{cases}
\ds |F(t,x,w)|&\leq~  \Gamma_1 \cdot \left[ \ell(t)\cdot \big(\big|x\ln|x| \big|+ t^{\alpha}\big) + t^{\alpha-1}\right], \\[3mm]
\ds |F_x(t,x,w)|&\leq~   \Gamma_1 \cdot \left[t^{\alpha-3/2}+\big(\ell(t) + t^{\alpha-1}\big)\cdot |x|^{\alpha-1}\right].
%\left[(|x|t)^{\alpha-1} + \ell(t)\left(t(|x|+t)^{\alpha-2} + \big|\ln|x|\big|\right)\right]\,.
\end{cases}
%\quad a.e.~t\in [0,T],
\eeq
Furthermore, for every $\delta>0$ sufficiently small, we obtain that 
\bel{eqF2}
\|F(t,\cdot,w)\|_{H^2(\R\backslash [-\delta,\delta])}~\leq~\Gamma_1\cdot \left[t^{-3/4}+  \left(t^{\alpha-1}+\ell(t)\right)\cdot\delta^{-3/4}\right], \ \ \text{with} \eeq
\bel{Gamma}
\Gamma_1~=~\O(1)\cdot\Gamma_f^2 C_0^2 (M_0^3+1) \left[1+b_0^{\alpha-3}\right],\ \Gamma_f~=~1+\|f\|_{\mathcal{C}^4([-2M_0-C_0,2M_0+C_0])}.
\eeq
\end{lemma}
{\bf Proof.} First, recall that we expressed $F$ in (\ref{vp-est}) in  terms of $A^{(w)}$, $B^{(w)}$ and $C^{(w)}$. 
\medskip

{\bf 1. Estimates on $B^{(w)}$.}  Applying Lemma \ref{B-est} for $w(t,\cdot)$, we get for every $0<|x|<1/4$ and $\delta>0$ that 
\[
\left|B^{(w)}(t,x)\right|~\leq~\O(1)\cdot M_0,\qquad  \left|B_x^{(w)}(t,x)\right|~\leq~\O(1)\cdot M_0\cdot\ln^2|x|,
\]
and 
\[
\left\|B^{(w)}(t,\cdot)\right\|_{H^2(\R\backslash [-\delta,\delta])}~\leq~\O(1)\cdot M_0\cdot \delta^{-2/3}.
\]

{\bf 2. Estimates on $A^{(w)}$.}  We are now giving some bounds on  $A^{(w)}$ by splitting it into two parts
\bel{A-split}
A^{(w)}~=~{\bf G}\left[\vp^{(w)}\right]  -A^{(w)}_1,\qquad A^{(w)}_1~=~d^{(w)}\cdot\vp^{(w)}_{x}
\eeq
with 
\[
d^{(w)}~\doteq~f'\left(w+\vp^{(w)}\right)-f'(w).
\]
From (\ref{sigma})-(\ref{ul-b}) and Lemma \ref{vp-est}, we estimate
\bel{d-es}
\begin{cases}
\left|\ds{d\over dx}d^{(w)}(t,x)\right|~\leq~ \ds\O(1) \cdot \Gamma_f C_0 \left(\big|\ln|x|\big| + \frac{1}{(|x|+t)^{1-\alpha}}\right),\\[4mm]
\left|\ds\frac{d^2}{dx^2} d^{(w)}(t,x)\right| ~ \leq ~\ds \O(1) \cdot \Gamma_f C_0 \left(\big|w_{xx}\big| + \frac{1}{|x|} + \frac{1}{(|x|+t)^{2-\alpha}}\right).
\end{cases}
\eeq
Notice that  $d^{(w)}(t,0)=0$, we derive from Lemma \ref{vp-est} that %let us write $\vp^{(w)}$ into two parts $\vp^{(w)}=\Tilde{\vp}^{(w)}+\ov{\vp}^{(w)}$ with 
\[
\begin{cases}
\left| A^{(w)}_1(t,x)\right|&~\leq~\ds \Gamma_f\cdot \big|\vp^{(w)}\vp^{(w)}_x\big|~\leq~\Gamma_fC_0^2\cdot |x|^{\alpha} \big(|t|^{\alpha-1}+\ln^2 |x|\big), \\[3mm]
\left|\ds{d\over dx} A^{(w)}_1(t,x)\right|&~\leq~\ds \O(1)\cdot \Gamma_f\cdot\left(\big|\vp^{(w)}\vp^{(w)}_{xx}\big|+\big|\vp^{(w)}_x\big|^2+\big|w_x\vp^{(w)}\vp^{(w)}_x\big|\right)\\
&~\leq~\ds\O(1) \cdot\Gamma_f C^2_0(M_0+1) \cdot |t x|^{\alpha-1}\,.
\end{cases}
\]
Moreover, keeping the leading order terms, one also gets  
\[
\begin{split}
\Big|\ds{d^2\over dx^2}& A^{(w)}_1(t,x)\Big|  \\
&\leq \O(1)\cdot \Gamma_f (M_0^2+1)\cdot \left[\big|\vp^{(w)}\vp_{xxx}^{(w)}\big|+\big|\vp_x^{(w)}\vp_{xx}^{(w)}\big|+\big|\vp_x^{(w)}\big|^3+\big|\vp^{(w)}\vp^{(w)}_{x}w_{xx}\big|\right]\\
&\leq\O(1)\cdot\Gamma_f (M_0^2+1)C^2_0\cdot \left[\left(\big|\ln|x|\big|+t^{\alpha-1}\right)\cdot\left({1\over |x|}+{1\over (|x|+t)^{2-\alpha}}\right)+{|w_{xx}|\over |x|^{1-\alpha}}\right].
\end{split}
\]
In particular, setting $\Gamma\doteq\Gamma_f^2 C_0^2 (M_0^3+1) \left[1+b_0^{\alpha-3}\right]$, we have 
\[
\begin{split}
\left\|A^{(w)}_1(t,\cdot)\right\|&_{H^2(\R\backslash [-\delta,\delta])}~\leq~\O(1)\Gamma\cdot \left(|\ln \delta|+t^{\alpha-1}\right) \left(\delta^{-1/2}+(\delta+t)^{\alpha-3/2}\right).
\end{split}
%\left({\delta^{-2/3}+(t+\delta)^{\alpha-3/2}\over t^{1-\alpha}}\right).
%\O(1)\Gamma_f (M_0^3+1)C^2_0\cdot \left({\delta^{-2/3}+(t+\delta)^{\alpha-3/2}\over t^{1-\alpha}}\right)\,.
\]
%that \textcolor{red}{(here)}
%\[
%\begin{split}
%\left\|\left[f'\left(w+\vp^{(w)}\right)-f'(w)\right]\cdot\vp^{(w)}_{x}(t,\cdot)\right\|_{H^2(\R\backslash [-\delta,\delta])}\\
%\leq~\O(1)\cdot \left({M^3_0\over b^3_0}+1\right)\cdot \left(\textcolor{red}{N_2}(\textcolor{red}{b_0}t)^{\alpha-1}\delta^{-2/3} + \textcolor{red}{N_3M_0}(\textcolor{red}{b_0}t+\delta)^{3\alpha-5/2}\right).
%\end{split}
%\]
To estimate the term ${\bf G}\left[\vp^{(w)}\right] (t,x)$ in $A^{(w)}$, we first recall  Lemma \ref{B-est} to obtain that
\[
\begin{cases}
\left|{\bf G}\big[\Tilde{\vp}^{(w)}(t,\cdot)\big](x)\right|\leq\O(1)\cdot \ds {M_0\over b_0},\qquad \left|\ds{d\over dx}{\bf G}\big[\Tilde{\vp}^{(w)}(t,\cdot)\big](x)\right|\leq\O(1)\cdot \ds {M_0\over b_0}\cdot\ln^2|x|,\\[4mm]
\left\|{\bf G}[\Tilde{\vp}^{(w)}(t,\cdot)]\right\|_{H^2(\R\backslash [-\delta,\delta])}\leq\O(1)\cdot\ds {M_0\over b_0}\cdot \delta^{-2/3}.
\end{cases}
\]
We observe that 
\[
\left|\overline{\vp}^{(w)}(t,0)\right| ~ = ~0, \ \ \text{and}\ \   \left|\frac{d^i}{dx^i} \overline{\vp}^{(w)}(t,x)\right| ~ \le ~ \O(1) \cdot (|x| + b_0 t)^{\alpha - i}, \; i \in \{1,2\},
\]
and therefore we obtain
\begin{equation*}
\left\|\overline{\vp}^{(w)}(t,\cdot)\right\|_{H^1(\R)} ~ \le ~ \O(1) \cdot \big(1 + b_0^{\alpha-1/2}\big) \cdot t^{\alpha-1/2},
\end{equation*} 
\begin{equation*}
\left\|\overline{\vp}^{(w)}(t,\cdot)\right\|_{H^2(\R)} ~ \le ~ \O(1) \cdot \big(1 + b_0^{\alpha-3/2}\big) \cdot t^{\alpha-3/2}.
\end{equation*}
Using the ${\bf L}^2$-continuity of ${\bf G}$, we then have that
\begin{align*}
\left\|{\bf G}\left[\overline{\vp}^{(w)}(t,\cdot)\right]\right\|_{L^\infty(\R)} \le 2 \cdot \left\|{\bf G}\left[\overline{\vp}^{(w)}(t,\cdot)\right]\right\|_{H^1(\R)} \le \O(1) \cdot \big(1 + b_0^{\alpha-1/2}\big) \cdot t^{\alpha-1/2},\\
\left\|\frac{d}{dx} {\bf G}\left[\overline{\vp}^{(w)}(t,\cdot)\right]\right\|_{L^\infty(\R)} \le 2 \cdot \left\|{\bf G}\left[\overline{\vp}^{(w)}(t,\cdot)\right]\right\|_{H^2(\R)} \le \O(1) \cdot \big(1 + b_0^{\alpha-3/2}\big) \cdot t^{\alpha-3/2}.
\end{align*}
Thus, for every $t>0$ and $|x|<1/4$, we have the following estimates:
\[\begin{split}
\left|A^{(w)}(t,x)\right|&~\leq~\O(1)\cdot\Gamma \cdot \left(|x|^{\alpha}\ln^2|x| +t^{\alpha-1}\right),\\
\left|{d\over dx}A^{(w)}(t,x)\right|&~\leq~\O(1)\cdot \Gamma \cdot t^{\alpha-1}\cdot \big(|x|^{\alpha-1}+t^{-1/2}\big),\\
\left\|A^{(w)}(t,\cdot)\right\|_{H^2(\R\backslash [-\delta,\delta])}&~\leq~\O(1)\cdot\Gamma\cdot t^{\alpha-1}\cdot\left(t^{-1/2}+\delta^{-2/3}+(\delta+t)^{\alpha-3/2}\right)\,.
\end{split}\]
%and
%\[
%\left\|A^{(w)}(t,\cdot)\right\|_{H^2(\R\backslash [-\delta,\delta])}~\leq~\O(1)\Gamma\cdot\left(t^{\alpha-1}\cdot (\delta+t)^{-3/2}\right)
%\]
%\[\begin{split}
%\left\|A^{(w)}(t,\cdot)\right\|_{H^2(\R\backslash [-\delta,\delta])}~\leq~\O(1)\cdot \Gamma_f (M_0^3+1) C_0^2 \cdot\\
%\left(t^{\alpha-1} \delta^{-2/3} + (t+\delta)^{\alpha-5/2} + \delta^{\alpha-1} + \delta^{-3/4}(\delta+t)^{\alpha-3/4} + t^{\alpha-1}(\delta+t)^{\alpha-3/2}\right).
%\end{split}\]
%medskip

{\bf 3. Estimates on $C^{(w)}$.}  From (\ref{F}) and (\ref{A})-(\ref{B}), we have that 
\[
C^{(\omega)}(t,x)~=~\vp^{(w)}_{t}(t,x)+b^{(w)}(t,x) \cdot \vp^{(w)}_{x}(t,x)+\sigma^{(\omega)}(t)\cdot\Lambda(x)\cdot \eta(x).
\]
%with
%\bel{btx}
%b(t,x)~:=~ f'\bigl(w(t,x)\bigr) - \frac{f\bigl(w^-(t)\bigr) - f\bigl(w^+(t)\bigr)}{w^-(t) - w^+(t)}.
%\eeq
%Observe that \textcolor{red}{(we don't need estimates on $b$ alone)}
%\bel{bbx}
%\begin{cases}
%|b(t,x)| + |b_x(t,x)| \le \O(1)\cdot \textcolor{red}{N_2} M_0\\
%|b_{xx}(t,x)| \le \O(1) \cdot \textcolor{red}{N_2}|w_{xx}(t,x)|
%\end{cases}
%\eeq
%and that
%\bel{bb}
%\left|b(t,x) - b_+^{(w)}(t)\right| \le \O(1) \cdot \textcolor{red}{N_2} \cdot M_0 \cdot x.
%\eeq
 Recalling (\ref{vp-tb}), we compute for $0<x<1/4$ that 
\[
\widetilde{\varphi}_x^{(w)}(t,x) =\ds {\sigma^{(w)}(t)\over b^{(w)}_{+}(t)}\cdot \left[\phi_x(x,0)-\phi_x\left(x,-t\cdot b^{(w)}_+(t)\right)\right],\quad \ov{\varphi}_x^{(w)}(t,x)=\overline{v}'\left(x-t\cdot b^{(w)}_+(t)\right)
\]
\[
\ov{\varphi}_t^{(w)}(t,x)=\left(b_+^{(w)}(t)+t\dot{b}_+^{(w)}(t)\right)\cdot\left[\overline{v}'\left(-t\cdot b^{(w)}_+(t)\right)-\overline{v}'\left(x-t\cdot b^{(w)}_+(t)\right)\right],
\]
\[
\begin{split}
\widetilde{\varphi}_t^{(w)}(t,x) &=\left({\dot{\sigma}^{(w)}(t)\over \sigma^{(w)}(t)}-{\dot{b}_+^{(w)}(t)\over b_+^{(w)}(t)}\right)\cdot \Tilde{\vp}^{(w)}(x,t) \\
&\qquad\qquad\qquad+ \frac{\sigma^{(w)}(t)}{b_+^{(w)}(t)}\cdot\left(b_+^{(w)}(t) + t \cdot \dot{b}_+^{(w)}(t)\right) \cdot \frac{d}{db}\phi\left(x,-tb_+^{(w)}(t)\right)
\end{split}
\]
%with 
%\bel{vp1}
%\Tilde{\varphi}_1(t,x)~:=~\left({\dot{\sigma}^{(w)}(t)\over \sigma^{(w)}(t)}-{\dot{b}_+^{(w)}(t)\over b_+^{(w)}(t)}\right)\cdot \Tilde{\vp}(x,t)
%\eeq
Since $\ds\frac{d}{db}\phi(x,b) = \Phi'(b)+\frac{d}{dx}\phi(x,b)$ and $\ds\frac{d}{dx}\phi(x,0)=-\Lambda(x)$, one has
\begin{multline}\label{C-1}
C^{(w)}(t,x)~=~\left({\dot{\sigma}^{(w)}(t)\over \sigma^{(w)}(t)}-{\dot{b}_+^{(w)}(t)\over b_+^{(w)}(t)}\right)\cdot \Tilde{\vp}^{(w)}(t,x) +\Big(b^{(w)}(t,x)-b^{(w)}_{+}(t)\Big)\cdot \varphi_x^{(w)}(t,x)\\
- t\dot{b}^{(w)}_+(t)\cdot \left[\ds {\sigma^{(w)}(t)\over b^{(w)}_{+}(t)}\cdot \Phi'\left(x-t\cdot b^{(w)}_+(t)\right)+\overline{v}'\left(x-t\cdot b^{(w)}_+(t)\right)\right]+E^{(w)}(t),
\end{multline}
with $E^{(w)}$ defined by 
\[
E^{(w)}(t)~\doteq~\left(b_+^{(w)}(t) + t \cdot \dot{b}_+^{(w)}(t)\right)\cdot \left[\frac{\sigma^{(w)}(t)}{b_+^{(w)}(t)}\cdot \Phi'\left(-t\cdot b^{(w)}_+(t)\right)+\overline{v}'\left(-t\cdot b^{(w)}_+(t)\right)\right].
\]
From (\ref{w-1-b})-(\ref{b-deri}) and Lemma \ref{vp-est} we obtain 
\[
\begin{split}
\left|C^{(w)}(t,x)\right|&~\le~ \O(1)\Gamma\cdot \big[\ell(t)\cdot \big(\big|x\ln|x|\big|+ t^{\alpha}\big)+t^{\alpha-1}\big]\,,\\
\left|C_x^{(w)}(t,x)\right|&~\le~ \O(1)\Gamma\cdot\left[\ell(t)\left[\frac{t}{(|x|+t)^{2-\alpha}} + \big|\ln |x|\big|\right] + \frac{1}{(|x|+t)^{1-\alpha}}\right]\,,\\
\left|C_{xx}^{(w)}(t,x)\right|&~ \le~ \O(1) \Gamma\cdot \left(\ell(t)\left[\frac{t}{(|x|+t)^{3-\alpha}} + \frac{1}{|x|}\right]  + (|w_{xx}|+1)\left[\big|\ln|x|\big| + \frac{1}{(|x|+t)^{1-\alpha}}\right]\right)\,,
\end{split}
\]
and 
\[
\left\|C^{(w)}(t,\cdot)\right\|_{H^2(\R\setminus[-\delta,\delta])} ~\le~ \O(1) \Gamma\cdot \left(t^{\alpha-1} + t\ell(t)(\delta+t)^{\alpha-5/2} + \ell(t)\delta^{-1/2}\right)\,.
\]
From the previous estimates on $A^{(w)}$, ${B}^{(w)}$, and (\ref{F1}), one then achieves (\ref{eqF1})-(\ref{eqF2}) for $\alpha>3/4$ . %The proof is complete.
%\[
%\left|C^{(w)}(t,x)\right| \le \O(1)\Gamma\cdot \left[\ell(t)+t^{\alpha-1}\right]
%\]
%\[
%\left|C_x^{(w)}(t,x)\right| \le \O(1) C_0 (M_0+1)\Gamma_f \left(\ell(t)\left[\frac{t}{(|x|+t)^{2-\alpha}} + \big|\ln |x|\big|\right] + \frac{1}{(|x|+t)^{1-\alpha}}\right)
%\]
%\[
%\left|C_{xx}^{(w)}(t,x)\right| \le \O(1) C_0 (M_0+1)\Gamma_f \left(\ell(t)\left[\frac{t}{(|x|+t)^{3-\alpha}} + \frac{1}{|x|}\right] + \frac{1}{(|x|+t)^{1-\alpha}} + |w_{xx}|\left[\big|\ln|x|\big| + \frac{1}{(|x|+t)^{1-\alpha}}\right]\right)
%\]
%\[
%\left\|C^{(w)}(t,\cdot)\right\|_{H^2(\R\setminus[-\delta,\delta])} ~\le~ \O(1) C_0 (M_0+1) \Gamma_f^2 b_0^{1-\alpha} \left(t^{\alpha-1} + t\ell(t)(\delta+t)^{\alpha-5/2} + \ell(t)\delta^{-1/2}\right)
%\]
\endproof
\medskip

To complete this section, we study the change in the function $F(t,x,w)$ as $w$ takes different values which   plays a key role in the proof of convergence of the approximations. More precisely, for any given two functions $w_1,w_2:[0,T]\times\R\to\R$ satisfying (\ref{w-1-b})-(\ref{Lip-w}),  we  provide a priori estimates on the difference 
\bel {F2-1}
{\bf F}^{(w_1,w_2)}(t,x)~\doteq~F(t,x,w_2)-F(t,x,w_1).
\eeq
in terms of $M_i(t)$ for $i\in\{1,2\}$ with
\bel{w2-1}
{\bf w}~\doteq~w_2-w_1,\quad\qquad M_{i}(t)~\doteq~\|{\bf w}(t,\cdot)\|_{H^i(\R\backslash\{0\})}.
\eeq
For all $x\in \R\backslash \{0\}$ and $t\in [0,T]$, it holds 
\bel{sigma-2-1}
|{\bf w}(t,x)|~\leq~2M_1(t),\qquad |{\bf w}_x(t,x)|~\leq~ 2M_2(t),
\eeq
For simplicity, recalling (\ref{speed1}) and  (\ref{vp-tb}), we set 
\bel{b-vp-21}
{\bf b}~\doteq~b^{(w_2)}-b^{(w_1)},\quad\qquad \Psi~\doteq~\vp^{(w_2)}-\vp^{(w_1)}.
\eeq
We first estimate ${\bf b}^{(w_1,w_2)}$ by direct computations
\bel{b-speed-2-1}
\begin{cases}
\big|{\bf b}(t,x)\big|~\leq~4\Gamma_f\cdot M_1(t),\qquad \big|{\bf b}_x(t,x)\big|~\leq~2 \Gamma_f\cdot \left(M_0M_1(t)+|{\bf w}_x(t,x)|\right),\\[4mm]
 \big|{\bf b}_{xx}(t,x)\big|~\leq~  4\Gamma_f\cdot \left(\left[M_0^2 +|w_{1,xx}(t,x)|\right]\cdot M_1(t)+2M_0M_2(t)+|{\bf w}_{xx}(t,x)|\right),
\end{cases}
\eeq
with $\Gamma_f$ given in (\ref{Gamma}). In  particular, this yields 
\bel{b-21-0}
\big|{\bf b}_{\pm}(t)\big|~\doteq~\big|b^{(w_2)}(t,0\pm)-b^{(w_2)}(t,0\pm)\big|~\leq~4 \Gamma_f\cdot M_1(t).
\eeq
Secondly, the difference of the corrector term  $ \Psi$ and its  transform ${\bf G}\big[ \Psi\big]$ is bounded by the following lemma.
\begin{lemma}\label{Psi21} For every $0<|x|<1/4$ and $0<t<\ds {1\over  4b_1}$, we have
\bel{Psi-2-1}
\begin{cases}
%{\color{blue}\left|{\bf G}\left[\Psi(t,\cdot)\right](x)\right|~\leq~ \Gamma_{f,1} M_1(t),}\qquad 
\left|\Psi(t,x)\right| \le  \Gamma_{f,1} M_1(t) |x|^{\alpha},\qquad \left|\ds{d\over dx}\Psi(t,x)\right|\leq\Gamma_{f,1} M_1(t)\left(\big|\ln|x|\big| + (|x|+t)^{\alpha-1}\right),\\[4mm]
\left|\ds{d^i\over dx^i}\Psi(t,x)\right|\leq\Gamma_{f,1} M_1(t) \left(|x|^{1-i} + (|x|+t)^{\alpha-i}\right),\quad i=2,3,\\[4mm]
\left\|{\bf G} \left[\Psi(t,\cdot)\right]\right\|_{H^1(\R\backslash\{0\})}, \|\Psi(t,\cdot)\|_{H^1(\R\backslash \{0\})}~\leq~  \Gamma_{f,1} M_1(t) t^{\alpha-1/2},
\end{cases}
\eeq
with $\Gamma_{f,1}=\O(1)(4b_1+16\Gamma_f M_0)(b_0^{\alpha-3}+b_0^{-2})$. Moreover, for every $\delta>0$ small, it holds 
\bel{Psi-22}
\begin{cases}
\left\|\Psi(t,\cdot)\right\|_{H^2(\R \setminus [-\delta,\delta])}~\leq~ \Gamma_{f,1}M_1(t)\left(t^{\alpha-3/2} + \delta^{-1/2}\right),\\[4mm]
\left\|{\bf G}\left[\Psi(t,\cdot)\right]\right\|_{H^2(\R \setminus [-\delta,\delta])}~\leq~ \Gamma_{f,1}M_1(t)\left(t^{\alpha-3/2} + \delta^{-2/3}\right).
\end{cases}
\eeq
\end{lemma}
{\bf Proof.} {\bf 1.} Recalling (\ref{vp-tb}), we first write $ \Psi=\ov{\Psi}+\Tilde{\Psi}$ with 
\bel{vp-21}
\ov{\Psi}~\doteq~\ov{\vp}^{(w_2)}-\ov{\vp}^{(w_1)},\quad\qquad \Tilde{\Psi}~\doteq~\Tilde{\vp}^{(w_2)}-\Tilde{\vp}^{(w_1)}.
\eeq
From (\ref{vp-tb}) and (\ref{Xadef}), for every $t\in [0,1/(4b_1)]$ and $|x|>0$, one has that $\overline{\Psi}(t,0)=0$, $\mathrm{supp}\left(\overline{\Psi}(t,\cdot)\right)\subseteq[-2,2]$, and 
%\[
%\overline{\Psi}^{\bf (w)}(t,0)~=~0,\qquad\mathrm{supp}\left(\overline{\Psi}^{\bf (w)}(t,\cdot)\right)~\subseteq~ [-2,2],
%\]
%and 
\[
\left|{d^i\over dx^i} \overline{\Psi}(t,x)\right|~\leq~\O(1)\cdot { M_1(t) \over (|x|+t b_0)^{i-\alpha}},\quad i=1,2,3,
\]
whence
\begin{equation*}
\begin{cases}
\left\|\overline{\Psi}(t,\cdot)\right\|_{H^1(\R)} ~ \le ~ \O(1)   \cdot \big(1+b_0^{\alpha-1/2}\big) \cdot t^{\alpha-1/2} \cdot M_1(t),\\[3mm]
\left\|\overline{\Psi}(t,\cdot)\right\|_{H^2(\R)} ~ \le ~ \O(1) \cdot \big(1+b_0^{\alpha-3/2}\big) \cdot t^{\alpha-3/2} \cdot M_1(t),
\end{cases}
\end{equation*}
and so the ${\bf L}^2$-continuity of ${\bf G}$ yields
\[
\left\|{\bf G}\left[\overline{\Psi}(t,\cdot)\right]\right\|_{{\bf L}^\infty(\R)} ~ \le ~ 2 \left\|{\bf G}\left[\overline{\Psi}(t,\cdot)\right]\right\|_{H^1(\R)} ~ \le ~ \O(1)  \cdot \big(1 + b_0^{\alpha-1/2}\big) \cdot t^{\alpha-1/2} \cdot M_1(t)
\]
and
\[
\left\|{\bf G}\left[\overline{\Psi}(t,\cdot)\right]\right\|_{H^2(\R)} ~ \le ~ \O(1)  \cdot \big(1 + b_0^{\alpha-3/2}\big) \cdot t^{\alpha-3/2} \cdot M_1(t).
\]
%\begin{align*}
%\left\|{\bf G}\left[\overline{\Psi}^{({\bf w})}(t,\cdot)\right]\right\|_{{L}^\infty(\R)} ~ \le ~ 2  \left\|{\bf G}\left[\overline{\Psi}^{({\bf w})}(t,\cdot)\right]\right\|_{H^1(\R)} ~ \le ~ \O(1)  \cdot (1 + b_0^{\alpha-1/2}) \cdot t^{\alpha-1/2} \cdot M_1(t),\\
%\left\|\frac{d}{dx} {\bf G}\left[\overline{\Psi}^{({\bf w})}(t,\cdot)\right]\right\|_{{L}^\infty(\R)} ~ \le ~ 2  \left\|{\bf G}\left[\overline{\Psi}^{({\bf w})}(t,\cdot)\right]\right\|_{H^2(\R)} ~ \le ~ \O(1) \cdot  (1 + b_0^{\alpha-3/2}) \cdot t^{\alpha-3/2} \cdot M_1(t).
%\end{align*}

{\bf 2.} Fix $0 < x < 1/4$ (similar computations hold for $-1/4 < x < 0$). To estimate $\Tilde{\Psi}$, we write
\begin{multline*}
\Tilde{\Psi}(t,x) = \left({\sigma^{(w_2)}(t) \over b^{(w_2)}_{+}(t)} - {\sigma^{(w_1)}(t) \over b^{(w_1)}_{+}(t)}\right) \biggl(\phi(x,0) - \phi\Bigl(x,-t\cdot b_+^{(w_1)}\Bigr)\biggr)\\
+ \frac{\sigma^{(w_2)}(t) {\bf b}_{+}(t)t}{b_+^{(w_2)}(t)} \cdot \int_{0}^{1} {d\over db}\phi\left(x,-t \cdot b_+^{(w_2)}(t) + \tau {\bf b}_{+}(t) t\right) d\tau.
\end{multline*}
From (\ref{ul-b}) and (\ref{b-21-0}), it holds
\bel{d-est}
\begin{split}
\left|{\sigma^{(w_2)}(t) \over b^{(w_2)}_{+}(t)} -{\sigma^{(w_1)}(t) \over b^{(w_1)}_{+}(t)}\right|&~\leq~{b_1\cdot |{\bf w}(t,0+)-{\bf w}(t,0-)|+4M_0\cdot \big|{\bf b}_{+}(t)\big|\over b_0^2}\\
&~\leq~ {4b_1+16\Gamma_f M_0\over b_0^2}\cdot M_1(t)~\doteq~\Gamma_{f,0}  M_1(t).
\end{split}
\eeq
Recalling (\ref{K1}), we  estimate 
\[
\begin{split}
\Big |\phi_x\left(x,-tb_+^{(w_1)}(t)\right) - &\phi_x(x,0)\Big| ~ = ~ \left|\int_0^{-tb_+^{(w_1)}(t)} K(x + \tau) \, d\tau\right| ~ \le ~ \O(1) \cdot \int_0^{b_1t} \frac{d\tau}{x + \tau}\\
&~=~\O(1)\cdot \ln \left(1+{b_1t\over x}\right)~\leq~\O(1)\cdot \min\left\{|\ln x|, \left(\frac{b_1 t}{x}\right)^{\alpha-1/2}\right\},
\end{split}
\]
and 
%\begin{align*}
%\left|\phi_x\left(x,-tb_+^{(w_1)}(t)\right) - \phi_x(x,0)\right| ~ = ~ \left|\int_0^{-tb_+^{(w_1)}} K(x + \tau) \, d\tau\right| ~ \le ~ \O(1) \cdot \int_0^{b_1t} \frac{d\tau}{x + \tau}\\
%= ~ \O(1) \cdot \int_0^{b_1t} \frac{d\tau}{(x+\tau)^{\alpha-1/2} (x+\tau)^{3/2-\alpha}} ~ \le ~ \O(1) \cdot \left(\frac{b_1 t}{x}\right)^{\alpha-1/2}.
%\end{align*}
%At the same time,
%\begin{equation*}
%\left|\phi_x\left(x,-tb_+^{(w_1)}(t)\right) - \phi_x(x,0)\right| ~ \le ~ \left|\phi_x\left(x,-tb_+^{(w_1)}(t)\right)\right| + \left|\phi_x(x,0)\right| ~ \le ~ \O(1) \cdot |\ln x|.
%\end{equation*}
%Moreover, since
\begin{align*}
\left|\frac{d}{dx}\int_{0}^{1} {d\over db}\phi\left(x,-t \cdot b_+^{(w_2)}(t) + \tau {\bf b}_{+}(t) t\right) d\tau\right| = \left|\int_0^1 K\left(x - t \cdot b_+^{(w_2)}(t) + \tau {\bf b}_{+}(t) t\right) d\tau\right|\\
\le \O(1) \cdot \int_0^1 \frac{d\tau}{x - (1-\tau) t b_+^{(w_2)}(t) - \tau t b_+^{(w_1)}} \le \O(1) \cdot \int_0^1 \frac{d\tau}{x + b_0 t} = \O(1) \cdot \frac{1}{x + b_0 t}.
\end{align*}
Thus, keeping into account that $\Tilde{\Psi}(t,0)=0$, we then have 
\begin{align*}
\left|\widetilde{\Psi}(t,x)\right| ~ \le ~ \O(1) \cdot \Gamma_{f,0} b_1^{\alpha-1/2} M_1(t) \cdot \min\big\{t^{\alpha-1/2} x^{3/2-\alpha},\big|x\ln x\big|\big\},\\
\left|\frac{d}{dx}\widetilde{\Psi}(t,x)\right| ~ \le \O(1) \cdot \Gamma_{f,0} b_1^{\alpha-1/2} b_0^{-1} M_1(t)\cdot \min\biggl\{\biggl(\frac{t}{x}\biggr)^{\alpha-1/2},\big|\ln x\big|\biggr\}.
\end{align*}
This, together with the ${\bf L}^2$-continuity of ${\bf G}$, yields
\begin{align*}
\left\|{\bf G}\left[\widetilde{\Psi}(t,\cdot)\right]\right\|_{{\bf L}^\infty(\R)} ~ & \le ~ 2\cdot \left\|{\bf G}\left[\widetilde{\Psi}(t,\cdot)\right]\right\|_{H^1(\R \setminus \{0\})} ~ \le ~ \O(1) \cdot \left\|\widetilde{\Psi}(t,\cdot)\right\|_{H^1(\R \setminus \{0\})}\\
& \le ~ \O(1) \cdot \Gamma_{f,0} b_1^{\alpha-1/2} M_1(t) t^{\alpha-1/2}.
\end{align*}
Concerning the other derivatives of $\widetilde{\Psi}$, direct computations, together with \eqref{ul-b}, \eqref{b-21-0}, and \eqref{d-est}, yield
\begin{equation*}
\left|\frac{d^i}{dx^i} \widetilde{\Psi}(t,x)\right| ~ \le ~ O(1) \cdot \Gamma_{f,1} M_1(t) \frac{1}{x^{i-1}}, \qquad i \in \{2,3\},
\end{equation*}
and, in particular,
\begin{equation*}
\left\|\widetilde{\Psi}(t,\cdot)\right\|_{H^2(\R \setminus [-\delta,\delta])} ~ \le ~ \O(1) \cdot \Gamma_{f,1} M_1(t) \delta^{-1/2}.
\end{equation*}

{\bf 3.} In order to estimate $\left\|{\bf G}\left[\Tilde{\Psi}(t,\cdot)\right]\right\|_{H^2(\R \setminus [-\delta,\delta])}$, we write $\Tilde{\Psi}=\Tilde{\Psi}_1+\Tilde{\Psi}_2$, with 
\begin{equation*}
\Tilde{\Psi}_1(t,x)=\left[\left({\sigma^{(w_2)}(t) \over b^{(w_2)}_{-}(t)} -{\sigma^{(w_1)}(t) \over b^{(w_1)}_{-}(t)}\right)\chi_{]-\infty,0[}+\left({\sigma^{(w_2)}(t) \over b^{(w_2)}_{+}(t)} -{\sigma^{(w_1)}(t) \over b^{(w_1)}_{+}(t)}\right) \cdot\chi_{]0,\infty[}\right] \phi(x,0)
\end{equation*}
From (\ref{phi-1}), (\ref{d-est}), and Lemma \ref{A1}, we obtain
\begin{align*}
%\left|\frac{d^i}{d x^i}\Tilde{\Psi}_1(t,x)\right|~\leq~\Gamma_{f,0}M_1(t)\cdot \left|\frac{d^i}{d x^i}\phi(x,0)\right|
\left\|\frac{d^2}{d x^2}{\bf G}[\Tilde{\Psi}_1(t,\cdot)]\right\|_{H^2(\R \setminus [-\delta,\delta])} ~ & \leq ~ \O(1) \cdot \Gamma_{f,0} M_1(t) \left\|\frac{d^2}{d x^2}{\bf G}[\phi(x,0)]\right\|_{H^2(\R \setminus [-\delta,\delta])}\\
& \le ~ \O(1) \cdot \Gamma_{f,0} M_1(t) \delta^{-2/3}.
\end{align*}
On the other hand, from (\ref{vp-tb}) and (\ref{vp-21}), for every $x>0$, we can write 
\begin{multline*}
\Tilde{\Psi}_2(t,x)~=~ \left({\sigma^{(w_1)}(t) \over b^{(w_1)}_{+}(t)} -{\sigma^{(w_2)}(t) \over b^{(w_2)}_{+}(t)}\right) \phi\left(x,-t \cdot b_+^{(w_1)}(t)\right)\\+\frac{\sigma^{(w_2)}(t) {\bf b}_{+}(t)t}{b_+^{(w_2)}(t)}\cdot\int_{0}^{1}{d\over db}\phi\left(x,-t \cdot b_+^{(w_2)}(t)+\tau {\bf b}_{+}(t) t\right)d\tau.
\end{multline*}
Combining \eqref{ul-b}, \eqref{b-21-0}, and \eqref{d-est}, we obtain 
\[
%\left|{d\over dx}\Tilde{\Psi}_2(t,x)\right|~\leq~ \O(1) \Gamma_{f,0} M_1(t) |\ln(|x|+b_0t)| ,\quad
\left|{d^2\over dx^2}\Tilde{\Psi}_2(t,x)\right|~\leq~  \O(1)\cdot { \Gamma_{f,0} M_1(t) \over (|x|+b_0t)}.
\] 
Similarly, the above estimates also hold for $x<0$. Thus, noticing that $\Tilde{\Psi}_2(t,0)=0$ and using again the ${\bf L}^2$-continuity of ${\bf G}$, we get 
\[
%\left\|\Tilde{\Psi}_2(t,\cdot)\right\|_{ H^1(\R)}~\leq~ \O(1) \Gamma_{f,0} M_1(t),
\left\|{\bf G}\left[\Tilde{\Psi}_2(t,\cdot)\right]\right\|_{H^2(\R)} ~ \le ~  \O(1)\cdot \left\|\Tilde{\Psi}_2(t,\cdot)\right\|_{ H^2(\R)}~\leq~ \O(1)\cdot \Gamma_{f,0} b_0^{-1/2} M_1(t) t^{-1/2}.
\]

Combining all the above estimates and using the fact that $3/4<\alpha<1$, we finally obtain (\ref{Psi-2-1}) and (\ref{Psi-22}).
%
%
%for all $x>0$. 
%\begin{multline*}
%\Tilde{\Psi}_2^{\bf (w)}(t,x)~=~ \left({\sigma^{(w_1)}(t) \over b^{(w_1)}_{+}(t)} -{\sigma^{(w_2)}(t) \over b^{(w_2)}_{+}(t)}\right) \phi\left(x,-t \cdot b_+^{(w_1)}(t)\right)\\
% + \frac{\sigma^{(w_2)}(t)}{b_+^{(w_2)}(t)} \left[\phi\left(x,-t \cdot b_+^{(w_1)}(t)\right) - \phi\left(x,-t \cdot b_+^{(w_2)}(t)\right)\right]
%\end{multline*}
%\begin{align*}
%\Tilde{\Psi}_2^{\bf (w)}(t,x) & = \left({\sigma^{(w_1)}(t) \over b^{(w_1)}_{\pm}(t)} -{\sigma^{(w_2)}(t) \over b^{(w_2)}_{\pm}(t)}\right) \phi\left(x,-t \cdot b_\pm^{(w_1)}\right)\\
%& \quad + \frac{\sigma^{(w_2)}(t)}{b_\pm^{(w_2)}(t)} \left[\phi\left(x,-t \cdot b_\pm^{(w_1)}(t)\right) - \phi\left(x,-t \cdot b_\pm^{(w_2)}(t)\right)\right]
%\end{align*}
\endproof
\medskip

Using the estimates in Lemma \ref{Psi21}, we provide an $H^2$ bound on ${\bf F}^{(w_1,w_2)}$ which allows us to obtain the convergence of a sequence of approximate solutions $w^{(k)}$ to (\ref{wne}) in Lemma \ref{w-n}.
\begin{lemma}\label{F2}
%Given two functions $w_1,w_2:[0,T]\times\R\to\R$ satisfying (\ref{w-1-b})-(\ref{Lip-w}),  set 
%\[
%{\bf w}~\doteq~w_2-w_1,\quad z^{\pm}(t)~\doteq~z(t,0\pm),\qquad M_{i}(t)~\doteq~\|z(t,\cdot)\|_{H^i(\R\backslash\{0\})}\quad i\in\{1,2\}.
%\
There exists a constant $\Gamma_2 \doteq\O(1)\cdot\Gamma_f \left(\Gamma_{f,1}^2 M_0^2 C_0^3  b_1 + 1\right)$ and $T>0$ sufficiently small such that for every $|x|<1/4$, $0<t<T$, and  $\delta>0$, we have 
\begin{equation}\label{F1-2}
\ds \left|{\bf F}^{(w_1,w_2)}(t,x)\right|\le \Gamma_2\cdot \Bigl(M_2(t) \cdot \left[t^{\alpha-1}+\ell(t)(t^\alpha + |x|^\alpha)\right] + \Sigma(t) |x|^{\alpha} \Bigr)
\end{equation}
\begin{equation}\label{F1-2-2}
\ds\left\|{\bf F}^{(w_1,w_2)}(t,\cdot)\right\|_{H^2(\R \setminus [-\delta,\delta])} \le \Gamma_2\cdot \Bigl( M_2(t)\cdot\left[\ell(t) \delta^{\alpha-3/2} + \delta^{-2/3} t^{2\alpha-11/6}\right]+ \Sigma(t) \delta^{\alpha-3/2}\Bigr),
\end{equation}
with $\Sigma(t) \doteq \max\left\{\left|\dot{w}_2^+(t) - \dot{w}_1^+(t)\right|,\left|\dot{w}_2^-(t) - \dot{w}_1^-(t)\right|\right\}$.
%\begin{equation*}
%\Sigma(t) \doteq \max\left\{\left|\dot{w}_2^+(t) - \dot{w}_1^+(t)\right|,\left|\dot{w}_2^-(t) - \dot{w}_1^-(t)\right|\right\}.
%\end{equation*}
%and
%\begin{equation*}
%\Gamma_2 ~ \doteq ~ \Gamma_f (\Gamma_{f,1}^2 M_0^2 C_0^3  b_1 + 1).
%\end{equation*}
%\bel{F1-2}
%\big|{\bf F}^{(w_1,w_2)}(t,x)\big|~\leq~????,\qquad \big\|{\bf F}^{(w_1,w_2)}(t,\cdot)\big\|_{H^2(\R\backslash [-\delta,\delta])}~\leq~?????
%\eeq
%Moreover, for every $\delta>0$, it holds 
%\bel{F-d-d} 
%\big\|{\bf F}^{(w_1,w_2)}(t,\cdot)\big\|_{H^2(\R\backslash [-\delta,\delta])}~\leq~?????
%\eeq
\end{lemma}
{\bf Proof.}  {\bf 1.} Recalling (\ref{F1}), (\ref{A-split}), and (\ref{w2-1})-(\ref{b-vp-21}),  we have 
\[
{\bf F}^{(w_1,w_2)}(t,x)~=~{\bf A}_2(t,x)-{\bf A}_1(t,x)+B^{({\bf w})}(t,x)-{\bf C}(t,x),
\]
with
\[
{\bf A}_2~\doteq~{\bf G}\big[\Psi\big],\qquad {\bf A}_1~\doteq~A_1^{(w_2)}-A_1^{(w_1)},\qquad  {\bf C}~\doteq~C^{(w_2)}-C^{(w_1)}.
\]
From Lemma \ref{B-est} and Lemma \ref{Psi21},  for every $0<|x|<1/4$ and $\delta>0$, it holds
\begin{equation*}
\begin{cases}
\left|B^{({\bf w})}(t,x)\right|~\leq~\O(1)\cdot M_1(t),\qquad \left\|B^{({\bf w})}(t,\cdot)\right\|_{H^{1}(\R\backslash\{0\})}~\leq~\O(1)\cdot M_1(t),\\[3mm]
\left\|B^{({\bf w})}(t,\cdot)\right\|_{H^2(\R\backslash [-\delta,\delta])}~\leq~\O(1)\cdot M_2(t)\cdot \delta^{-2/3},
\end{cases}
\end{equation*}
and 
\[
\begin{cases}
\left|{\bf A}_2(t,x)\right|~\leq~\O(1) \cdot \Gamma_{f,0} M_1(t),\qquad \left\|{\bf A}_2(t,\cdot)\right\|_{H^{1}(\R\backslash\{0\})}~\leq~\O(1)\cdot\Gamma_{f,1} M_1(t),\\[3mm]
\left\|{\bf A}_2(t,\cdot)\right\|_{H^2(\R\backslash [-\delta,\delta])}~\leq~\O(1) \cdot \Gamma_{f,1}M_1(t)\left(t^{\alpha-3/2} + \delta^{-2/3}\right).
\end{cases}
\]

{\bf 2.} To bound the term $\bf A_1$, we first recall  $d^{(w)}\doteq f'\left(w+\vp^{(w)}\right)-f'(w)$ and write 
\[
{\bf A}_1~=~d^{(w_1)}\Psi_x+{\bf d}\cdot \vp_x^{(w_2)}~\doteq~ {\bf A}_{1,1}+{\bf A}_{1,2},\qquad {\bf d}~\doteq~d^{(w_2)}-d^{(w_1)}~=~{\bf d}_1+{\bf d}_2\,,
\]
with 
\[
\begin{cases}
{\bf d}_1~=~\ds\left(\int_{0}^{1}f''\left(w_1+\tau  \cdot \vp^{(w_1)}\right) d\tau\right)\cdot\Psi,\\[4mm]
{\bf d}_2~=~\ds\vp^{(w_2)}\cdot \int_{0}^{1}\left[\int_{0}^{1}f^{(3)}\left(w_1+\tau  \cdot \vp^{(w_1)}+[{\bf w}+\tau\Psi]\cdot s\right)ds\cdot \left({\bf w}+\tau \Psi\right)\right]d\tau.
\end{cases}
\]
Since $\|w_i\|_{H^2(\R\backslash \{0\})}\leq M_0$ and $f\in {C}^4$, we directly estimate
\bel{d-w}
%\begin{cases}
\left|\ds{d\over dx}{\bf d}^{({\bf w})}(t,x)\right|\leq \mathcal{C}\cdot\min\left\{M_1(t) |x|^{\alpha-1},M_2(t) \left(\big|\ln|x|\big| + \frac{1}{(|x|+t)^{1-\alpha}}\right)\right\}
\eeq
\bel{d-w2}
%\left|\ds{d\over dx}{\bf d}^{({\bf w})}(t,x)\right|~\leq~\O(1) \cdot \ds\Gamma_f \Gamma_{f,1} C_0^2  M_2(t) \left(\big|\ln|x|\big| + \frac{1}{(|x|+t)^{1-\alpha}}\right),\\[4mm]
\ds\left|\frac{d^2}{dx^2} {\bf d}^{\bf (w)}(t,x)\right|\leq \mathcal{C} \cdot \ds\left[|{\bf w}_{xx}|+ M_2(t)\left(\big|w_{2,xx}\big|+{1\over |x|}+{1\over (|x|+t)^{2-\alpha}}\right)\right],
%\end{cases}
\eeq
where $\mathcal{C} = \O(1) \cdot \Gamma_f \Gamma_{f,1} C_0^2$. \\

Notice that  $d^{(w_1)}(t,0)={\bf d}^{({\bf w})}(t,0)=0$. By keeping the leading order terms, we derive from (\ref{d-es}), Lemma \ref{vp-est} and Lemma \ref{Psi21} that for every $|x|<1/4$ and $0<t<T$, it holds
\[
%\begin{cases}
\bigg|{\bf A}_{1}(t,x)\bigg|~ \leq~ \left|{\bf A}_{1,1}(t,x)\right|+\left|{\bf A}_{1,2}(t,x)\right|~\leq~\ \Gamma_2M_2(t)\cdot |x|^{\alpha}t^{\alpha-1},
\]
\[
\bigg|\ds{d\over dx}{\bf A}_{1}(t,x)\bigg|~ \leq~\left|{\bf A}_{1,1}(t,x)\right|+\left|{\bf A}_{1,2}(t,x)\right|~\leq~  \Gamma_2 M_2(t)\cdot \left(\ln^2|x|+{t^{\alpha-1}\over (|x|+t)^{1-\alpha}}\right),
%\end{cases}
\]
and 
\begin{align*}
\bigg|\ds{d^2\over dx^2}{\bf A}_{1}(t,x)\bigg|~\leq~ \Gamma_2\cdot \biggl[M_2(t) \left(\big|\ln|x|\big| + t^{\alpha-1}\right) \left(\frac{1}{|x|} + \frac{1}{(|x|+t)^{2-\alpha}}\right)\\
+ \biggl(|{\bf w}_{xx}| + M_2(t)(|w_{1,xx}|+|w_{2,xx}|)\biggr) \left(\big|\ln|x|\big| + \frac{1}{(|x|+t)^{1-\alpha}}\right)\biggr],
\end{align*}
and thus
\begin{equation*}
\left\|{\bf A}_1(t,\cdot)\right\|_{H^2(\R \setminus [-\delta,\delta])} ~ \le ~ \Gamma_2 M_2(t) \Bigl[(\ln\delta + t^{\alpha-1}) \bigl(\delta^{-1/2} + (\delta+t)^{\alpha-3/2}\bigr)\Bigr].
\end{equation*}

{\bf 3.} It remains to estimate ${\bf C}$. From (\ref{C-1}), (\ref{w2-1}), (\ref{sigma-2-1}), (\ref{b-21-0}), we have 
\begin{multline*}
\big|{\bf C}(t,0+)\big|~\leq~\Big|\sigma^{(w_2)}\cdot \Phi'\big(-tb_+^{(w_2)}\big)-\sigma^{(w_1)}\cdot \Phi'\big(-tb_+^{(w_1)}\big)\Big|\\
+\Big|b_+^{(w_2)}\cdot \bar{v}'\big(-tb_{+}^{(w_2)}\big)-b_+^{(w_1)}\cdot \bar{v}'\big(-tb_{+}^{(w_1)}\big)\Big|~\leq~ \Gamma_2 M_1(t) t^{\alpha-1}.
\end{multline*}
Moreover, for  $0<x<1/4$ (similar estimates hold when $-1/4<x<0$), it holds
\[
{d \over dx}{\bf C}(t,x) ~=~ {d \over dx}{\bf C}_1(t,x)+{d \over dx}{\bf C}_2(t,x)+{d \over dx} {\bf C}_3(t,x),
\]
with 
\[
%\begin{cases}
{\bf C}_1=\left[\frac{\dot{\sigma}^{(w_2)}(t)}{\sigma^{(w_2)}(t)} - \frac{\dot{\sigma}^{(w_1)}(t)}{\sigma^{(w_1)}(t)} -\frac{\dot{b}_+^{(w_2)}(t)}{b_+^{(w_2)}(t)} + \frac{\dot{b}_+^{(w_1)}(t)}{b_+^{(w_1)}(t)}\right]\widetilde \vp^{(w_2)}+\ds \left(\frac{\dot{\sigma}^{(w_1)}(t)}{\sigma^{(w_1)}(t)} - \frac{\dot{b}_+^{(w_1)}(t)}{b_+^{(w_1)}(t)}\right) \widetilde\Psi,
\]
\[
{\bf C}_2 = \left({\bf b}(t,x)-{\bf b}_+(t)\right)\cdot \vp_x^{(w_2)}+\left(b^{(w_1)}(t,x)-b^{(w_1)}_+(t)\right)\cdot \Psi_x,
\]
\[
{\bf C_3}=-\ds t \dot{{\bf b}}_+(t)\cdot D^{(w_2)}+ t\dot{b}_+^{(w_1)}(t)\cdot \left[D^{(w_1)}-D^{(w_2)}\right],
%\end{cases}
\]
and 
\[
D^{(w)}(t,x)~\doteq~\frac{\sigma^{(w)}(t)}{b_+^{(w)}(t)} \Phi'\left(x-tb_+^{(w)}(t)\right) + \overline{v}'\left(x-tb_+^{(w)}(t)\right).
\]
%
%\begin{align*}
%{\bf C}_1(t,x) & = \left(\frac{\dot{\sigma}^{(w_2)}(t)}{\sigma^{(w_2)}(t)} - \frac{\dot{b}^{(w_2)}(t)}{b^{(w_2)}(t)}\right) \widetilde{\varphi}_x^{(w_2)}(t,x) - \left(\frac{\dot{\sigma}^{(w_1)}(t)}{\sigma^{(w_1)}(t)} - \frac{\dot{b}^{(w_1)}(t)}{b^{(w_1)}(t)}\right) \widetilde{\varphi}_x^{(w_1)}(t,x)\\
%{\bf C}_2(t,x) & = \left(b^{(w_2)}(t,x) - b_\pm^{(w_2)}(t)\right) \vp_x^{(w_2)}(t,x) - \left(b^{(w_1)}(t,x) - b_\pm^{(w_1)}(t)\right) \vp_x^{(w_1)}(t,x)\\
%{\bf C}_3(t,x) & = - t\dot{b}_\pm^{(w_2)} \left[\frac{\sigma^{(w_2)}(t)}{b_\pm^{(w_2)}(t)} \Phi'\left(x-tb_\pm^{(w_2)}(t)\right) + \overline{v}'\left(x-tb_\pm^{(w_2)}(t)\right)\right]\\
%& \quad + t\dot{b}_\pm^{(w_1)} \left[\frac{\sigma^{(w_1)}(t)}{b_\pm^{(w_1)}(t)} \Phi'\left(x-tb_\pm^{(w_1)}(t)\right) + \overline{v}'\left(x-tb_\pm^{(w_1)}(t)\right)\right]\\
%{\bf C}_4(t,x) & = {\bf C}_4(t) = E^{(w_2)}(t) - E^{(w_1)}(t),
%\end{align*}
%where $E^{(w)}(t)$ is defined in \eqref{E}. Next, writing
%\begin{align*}
%{\bf C}_1(t,x) = \left[\left(\frac{\dot{\sigma}^{(w_2)}(t)}{\sigma^{(w_2)}(t)} - \frac{\dot{\sigma}^{(w_1)}(t)}{\sigma^{(w_1)}(t)}\right) - \left(\frac{\dot{b}_\pm^{(w_2)}(t)}{b_\pm^{(w_2)}(t)} - \frac{\dot{b}_\pm^{(w_1)}(t)}{b_\pm^{(w_1)}(t)}\right)\right] \vp^{(w_2)}(t,x)\\
%+ \left(\frac{\dot{\sigma}^{(w_1)}(t)}{\sigma^{(w_1)}(t)} - \frac{\dot{b}_\pm^{(w_1)}(t)}{b_\pm^{(w_1)}(t)}\right) \Psi(t,x),
%\end{align*}
Recalling (\ref{w-1-b}), (\ref{ul-b}), (\ref{sigma-2-1}), and (\ref{b-21-0}), and taking $\tilde C = \O(1)  \frac{\Gamma_f(M_0 + b_1) + 1}{b_0^2}$,  we have
\[
\left|\frac{\dot{\sigma}^{(w_2)}(t)}{\sigma^{(w_2)}(t)} - \frac{\dot{\sigma}^{(w_1)}(t)}{\sigma^{(w_1)}(t)}\right| +\left|\frac{\dot{b}_+^{(w_2)}(t)}{b_+^{(w_2)}(t)} - \frac{\dot{b}_+^{(w_1)}(t)}{b_+^{(w_1)}(t)}\right|~\leq~  \tilde C \Bigl(\ell(t) M_1(t) + \Sigma(t)\Bigr),
\]
and this implies 
\[
\left|{d^i\over dx^i } {\bf C}_1\right|~\leq~ \Gamma_2\cdot \left[\Bigl(\ell(t) M_1(t) + \Sigma(t)\Bigr)\cdot \left|{d^i\over dx^i } \widetilde\vp^{(w_2)}\right|+ \Gamma \cdot \ell(t) \cdot \left|{d^i\over dx^i } \widetilde\Psi\right|\right].
\]
Moreover, from (\ref{b-speed}) and (\ref{b-speed-2-1}), taking $\hat C = \O(1) (\Gamma_f M_0^2 + \Gamma_f + b_1 + 1)$,  we obtain 
\[
\begin{cases}
\ds\left|{d\over dx}{\bf C}_2\right|~\leq~ \hat C \sum_{i=1}^{2} x^{i-1} \left(M_2(t) \left|{d^i\over dx^i}\vp^{(w_2)}\right|+\left|{d^i \over dx^i} \Psi\right|\right),\\[5mm]
\ds\Big|{d^2\over dx^2}{\bf C}_2\Big|~\leq~ \hat C \bigg[\sum_{i=2}^{3} x^{i-2} \left(M_2(t) \left|{d^i\over dx^i}\vp^{(w_2)}\right|+\left|{d^i\over dx^i}\Psi\right|\right)\\
\qquad\qquad +\bigl[|{\bf w}_{xx}| + M_2(t)(|w_{2,xx}|+1)\bigr] \left|\vp_x^{(w_2)}(t,x)\right| + \bigl(|w_{2,xx}|+1\bigr) \left|\Psi_x\right|\bigg].
\end{cases}
\]

In particular,  Lemma \ref{vp-est} and Lemma \ref{Psi21} imply
\bel{C12}
\begin{cases}
\ds \left|{d\over dx } {\bf C}_1\right|\leq \Gamma_2 \cdot \bigl[\ell(t) M_1(t) + \Sigma(t)\bigr] |\ln x|,\ \ds \left|{d^2\over dx^2 } {\bf C}_1\right|\leq \Gamma_2 \cdot {\ell(t) M_1(t) + \Sigma(t)\over x}, \\[4mm]
\ds \left|{d\over dx } {\bf C}_2\right|\leq \Gamma_2  M_2(t)\cdot \left(|\ln x| +(x+t)^{\alpha-1}\right),\\[4mm]
\ds \left|{d^2\over dx^2 } {\bf C}_2\right|\leq \Gamma_2\cdot \Big(M_2(t)\left(x^{-1}+(x+t)^{\alpha-2}\right) + \left(\left|{\bf w}_{xx}\right| + M_2(t)\left|w_{2,xx}\right|\right)\tilde C(x,t)\Big),
\end{cases}
\eeq
with $\tilde C(x,t) = |\ln x| +(x+t)^{\alpha-1}$. 
Noticing that
\[
\big|\dot{{\bf b}}_{+}(t)\big|~\leq~ \O(1)  \Gamma_f (M_0+1)  \bigl[\ell(t) M_1(t) + \Sigma(t)\bigr],\qquad \left|\dot{b}_+^{(w_1)}(t)\right|~\leq~ \O(1)  \Gamma_f (M_0+1) \cdot \ell(t),
\] 
we have 
\[
\left|{d^i\over dx^i } {\bf C}_3\right|\leq \tilde C t\left(\bigl[\ell(t) M_1(t) + \Sigma(t)\bigr] \left|{d^i\over dx^i } D^{(w_2)}\right|+ \ell(t) \left|{d^i\over dx^i } \left(D^{(w_2)}-D^{(w_1)}\right)\right|\right).
\]
with $\tilde C = \O(1)  \Gamma_f (M_0+1)$.  Thus, a direct computation yields
\[
\begin{cases}
\ds\left|{d^i\over dx^i } D^{(w_2)}(t,x)\right|~\leq~\O(1) \cdot \left(1+b_0^{\alpha-1-i}\right)\frac{1}{(x+t)^{1-\alpha+i}},\\[4mm]
\ds \left|{d^i\over dx^i } \left(D^{(w_2)}(t,x)-D^{(w_1)}(t,x)\right)\right|~\leq~ \O(1) \cdot \Gamma_f \left(1+b_0^{\alpha-2-i}\right)  \frac{M_1(t)}{(x+t)^{1-\alpha+i}},
\end{cases}
\]
and
\[
\ds \left|{d^i\over dx^i } {\bf C}_3(t,x)\right|~\leq~\Gamma_2 \cdot \frac{\ell(t) M_1(t) + \Sigma(t)}{(x+t)^{i-\alpha}}.
\]
Finally, combining the above estimates and (\ref{C12}), we obtain
\[
\begin{cases}
\ds\ds \left|{d\over dx } {\bf C}(t,x)\right|~\leq~ \Gamma_2 \cdot \Bigl[\ell(t) M_1(t) + M_2(t) + \Sigma(t)\Bigr] \left(|\ln x| + \frac{1}{(x+t)^{1-\alpha}}\right),\\[4mm]
\ds \left|{d^2\over dx^2 } {\bf C}(t,x)\right|~\leq~ \Gamma_2 \cdot \Bigl\{\Bigl[\ell(t) M_1(t) + M_2(t) + \Sigma(t)\Bigr] \left(\frac1x + \frac{1}{(x+t)^{2-\alpha}}\right)\\
\qquad \qquad \qquad\qquad \qquad \qquad \ds + \Bigl(|{\bf w}_{xx} + M_2(t)|w_{2,xx}|\Bigr) \left(|\ln x| + \frac{1}{(x+t)^{1-\alpha}}\right)\Bigr\},\\[4mm]
\left\|{\bf C}(t,\cdot)\right\|_{H^2(\R \setminus [-\delta,\delta])}~\leq~ \Gamma_2 \cdot \Bigl[\ell(t) M_1(t) + M_2(t) + \Sigma(t)\Bigr] \left(\delta^{-1/2} + (\delta+t)^{\alpha-3/2}\right),
\end{cases}
\]
and then we obtain the desired estimates.
\endproof

\section{Local existence and uniqueness of an entropic solution}
\setcounter{equation}{0}
Throughout this section, we give a proof of    Theorem~\ref{t:main} by constructing a local solution to the Cauchy problem \eqref{BL-1} with general initial data of the form $\ov u=\ov{w}+\ov v$ as in \eqref{IC1} where $\ov v$ is in  $\X_\alpha$ defined in (\ref{Xadef}) and $\ov{w}\in H^2({\R\backslash\{0\}})$ satisfies 
\bel{cond}
\delta_0~\doteq~\overline{w}(0-)-\overline{w}(0+)~>~0,\qquad {M_0\over 2}~\doteq~\|\overline{w}\|_{H^2(\R\backslash \{0\})}~<~\infty.
\eeq
The solution will be obtained as limit of a Cauchy sequence of approximate solutions $w_n(t,x)$ in $L^{\infty}([0,T],H^1(\R\backslash\{0\}))$ for some $T>0$ sufficiently small, following the two steps (i)-(ii) outlined at the end of Section~\ref{sec:main}. Indeed, we  first establish the existence and uniqueness of solutions to the Cauchy  problem (\ref{wne}).
\subsection{Construction of approximate solutions}
For fixed $n\geq 1$, let $w_n:[0,T]\times\R\to \R$ be a given function such that 
\bel{wn}
\|w_n(t,\cdot)\|_{H^2(\R\backslash\{0\})}~\leq~M_0,\qquad \left|w_n(t,0\pm)-\ov w(0\pm)\right|~\leq~{\delta_0\over 3}.
\eeq
For simplicity, recalling the definition of the corrector function $\vp^{(w)}(t,x)$ in (\ref{vp1}), (\ref{a}), and (\ref{F}),   we set
\bel{nota-n}
\begin{cases}
\sigma_n(t)~\doteq~w^-_n(t)-w^+_n(t),\qquad w_n^{\pm}(t)~\doteq~w_n(t,0\pm)\,,\\[4mm]
\vp_n(t,x)~\doteq~\vp^{(w_n)}(t,x),\qquad a_n(t,x)~\doteq~a(t,x,w_n).
\end{cases}
\eeq
The Cauchy problem (\ref{wne}) can be rewritten as
\bel{BH-3}
w_t+a_n(t,x)\cdot w_x~=~F(t,x,w),\qquad w(0,x)~=~\ov w(x).
\eeq 
Hence, to complete step (i), we shall prove that (\ref{BH-3}) admits a unique solution $w_{n+1}$ which  satisfies the bounds listed in (\ref{wn}). The construction of $w_{n+1}$ is divided into three steps:
\v

{\bf Step 1.} Let $t\mapsto x_n(t;\tau,x_\tau)$ be the solution to the Cauchy problem
\[
\dot{x}(t)~=~a_n(t,x(t)),\qquad x(\tau)~=~x_\tau,
\]
with $a_n(t,x)$ being the characteristic speed of (\ref{BH-3}), i.e.
\[
a_n(t,x)~=~b^{(w_n)}(t,x)+f'\left(w_n+\vp_n\right)-f'(w_n).
\]
The solution $w_{n+1}$ will be constructed by considering a sequence of approximate solutions $w^{(k)}$ to (\ref{wne}),
inductively defined as follows. 
\begin{enumerate}
\item  $w^{(1)}(t,\cdot)~\doteq~\overline{w}(\cdot)$ for all $t\geq 0$.
\item For every $k\geq 1$, $w^{(k+1)}(t,\cdot)$ solves the linear equation 
\bel{ln-wk}
w_t+a_n(t,x)\cdot w_x~=~F^{(k)}(t,x)~\doteq~F\left(t,x,w^{(k)}\right),\qquad w(0,\cdot)~=~\overline{w}(\cdot).
\eeq
 Equivalently, $w^{(k+1)}$
satisfies the integral identity
\bel{w-k+1}
w^{(k+1)}(t_0,x_0)~=~\overline{w}(x_n(0;t_0,x_0))+\int_{0}^{t_0}F^{(k)}(t,x_n(t;t_0,x_0))dt.
\eeq
\end{enumerate}
Recall the definition of $b_0$ and $b_1$ in (\ref{b0def}) and (\ref{b-speed}), we denote by 
\bel{I}
I^{\tau}_t~\doteq~[-b_0(\tau-t)/2,b_0(\tau-t)/2]\,,\qquad 0\leq t\leq \tau<\infty.
\eeq
%
%For every $0\leq t\leq \tau<\infty$, we denote by 
%\bel{I}
%I^{\tau}_t~\doteq~[-b_0(\tau-t)/2,b_0(\tau-t)/2]\,,
%\eeq
%with $b_0$ being defined in (\ref{b0def}).
%{\color{red}where $b_0$ is defined in (\ref{b0def}).}
The next lemma provides some properties including the Lipschitz continuous dependence of the characteristic curves $t\mapsto x_n(t;\tau,x_\tau)$. %{\color{red}Let us also recall from (\ref{b-speed}) the definition of $b_1$}
\begin{lemma}\label{Cr} Assume that $w_n$ and $\vp_n$ satisfy (\ref{wn})-(\ref{nota-n}). Then there exist   constants $\delta_1, T, K>0$ depending only on $M_0,\delta_0$ and $f$ such that for all 
$(x_\tau,\tau)\in \left([-\delta_1,\delta_1]\backslash\{0\}\right)\times (0,T]$, and for every $0\leq t\leq \tau$, we have 
\bel{x-n}
x_n(t;\tau,x_\tau)~\notin~I^{\tau}_t,\ \ \text{and}\ \  {b_0(\tau-t)\over 2}~\leq~\big|x_n(t;\tau,x_\tau)-x_\tau\big|~\leq~2b_1(\tau-t).
\eeq
Moreover, for any $0<x_1<x_2\leq \delta_1$ or $-\delta_1\leq x_1<x_2<0$, one has
\bel{L-x-r}
\big|x_n(t;\tau,x_1)-x_n(t;\tau,x_2)\big|~\leq~K\cdot |x_1-x_2|.
\eeq
\end{lemma}
{\bf Proof.} From Lemma \ref{vp-est}, it holds  for every $0<|x|<1/4$ and  $0<t<\ds {1\over  4b_1}$ that
\bel{dx-a}
\begin{split}
\left|{d\over dx}a_n(t,x)\right|&~\leq~\left|b_x^{(w_n)}(t,x)\right|+\Gamma_f\cdot \left(\big|\vp_{n,x}(t,x)\big|+2\cdot\big|w_{n,x}(t,x)\big|\right)\\
&~\leq~b_1+\Gamma_f\cdot \left(C_0\cdot |x|^{\alpha-1}+4M_0\right).% \le \Gamma_f C_0 |x|^{\alpha-1}.
\end{split}
\eeq
Recalling (\ref{ul-b}) and since $\varphi^{(w)}(t,0) = 0$, we have that 
\[
b_0~\leq~a_n(t,0-)~=~b_{-}^{(w_n)}(t)~\leq~b_1,\qquad -b_1~\leq~a_n(t,0+)~=~b^{(w_n)}_+(t)~\leq~-b_0,
\]
there exists a constant  $\delta_1>0$ depending on $b_1, b_0, C_0, M_0$, and $\Gamma_f$  such that 
\bel{bd-an}
\begin{cases}
-2b_1~\leq~a_n(t,x)~\leq~\ds-{b_0\over 2},&~~~ (t,x)\in [0,1/(4b_1)]\times \,]0,2\delta_1],\\[3mm]
\ds{b_0\over 2}~\leq~a_n(t,x)~\leq~2b_1,&~~~(t,x)\in [0,1/(4b_1)]\times [-2\delta_1,0[\,.
\end{cases}
\eeq
In particular, set $T\doteq \min\{1/(4b_1),\delta_1/(2b_1)\}$.  For $(x_\tau,\tau)\in \left([-\delta_1,\delta_1]\backslash\{0\}\right)\times (0,T]$, one has
\[
\big|x_n(t;\tau,x_\tau)\big|~\leq~ 2b_1|\tau|+|x_{\tau}|~\leq~2\delta_1,\qquad 0\leq t\leq\tau,
\]
and (\ref{bd-an}) implies (\ref{x-n}).
\medskip

To complete the proof, we shall establish (\ref{L-x-r}) for $0<x_1<x_2\leq \delta_1$, the other case being entirely similar. Set $z(t)\doteq\big|x_n(t;\tau,x_1)-x_n(t;\tau,x_2)\big|$. From (\ref{dx-a}), one has 
\[
\begin{split}
\left|{d\over dt}z(t)\right|&~\leq~\left(b_1+\Gamma_f\cdot \left[C_0\cdot |x_n(t;\tau,x_2)|^{\alpha-1}+4M_0\right]\right)\cdot z(t)\\
&\leq~\left(b_1+\Gamma_f\cdot \left[C_0\cdot \left|{b_0(\tau-t)\over 2}\right|^{\alpha-1}+4M_0\right]\right)\cdot z(t),
\end{split}
\]
%\[
%\left|{d\over dt}z(t)\right|~\leq~ \Gamma_f \cdot C_0 \cdot \left|x_n(t;\tau,x_2)\right|^{\alpha-1} \cdot z(t) ~ \le ~ \Gamma_f \cdot C_0 \cdot \biggl(\frac{b_0(\tau-t)}{2}\biggr)^{\alpha-1} \cdot z(t)
%\]
and this yields (\ref{L-x-r}).
\endproof
\medskip

As a consequence, for every $\tau\in [0,T]$, all characteristics starting at time $t=0$ inside the interval $I^{\tau}_0$ hit the origin before time $\tau$. On the other hand, since $a_n=b^{(w_n)}+d^{(w_n)}$, from (\ref{b-speed}), (\ref{d-es}), and (\ref{I}),  there exists a constant  $\Tilde{C}_0>0$ depending only on $M_0,\delta_0$, $\alpha$ and $f$ such that for every $x\in\R\backslash I^{\tau}_{t}$ and $t\in [0,\tau]$, it holds 
\bel{bound-a}
\big|a_{n,x}(t,x)\big|~\leq~\Tilde{C}_0\cdot (\tau-t)^{\alpha-1},\qquad \big|a_{n,xx}(t,x)\big|~\leq~\Tilde{C}_0\cdot \left(|w_{n,xx}(t,x)|+ |x|^{\alpha-2}\right).
\eeq

Hence, one can use  the same arguments as in \cite[Lemma 4.1]{BZ} to prove the following Lemma:
\begin{lemma}\label{est-l1} Under the same assumptions as in Lemma \ref{Cr}, there  exists  a constant $ T>0$ depending only on $M_0,\delta_0$ and $f$ such that for every $\tau\in [0,T]$ and any solution $v$ of the linear equation 
\bel{v-lem42}
v_t+a_n(t,x)\cdot v_x~=~0,\qquad v(0,\cdot)~=~\bar{v}~\in~ H^2\bigl(
\R\backslash I^{\tau}_0\bigr),
\eeq
one has
\bel{32-ineq}
\|v(\tau,\cdot)\|_{H^2\bigl(\R\backslash\{0\}\bigr)}~\leq~{3\over 2}\cdot \|\bar{v}\|_{H^2
(\R\backslash I^{\tau}_0)}.
\eeq
\end{lemma}
{\bf Proof.} To establish (\ref{32-ineq}) for a fixed $\tau\in [0,T]$, we set 
\[
Z(t)~\doteq~\|v(t,\cdot)\|^2_{H^2(\R\backslash I^{\tau}_t)}\qquad\forall t\in [0,\tau].
\]
%Notice that    $a_n=b^{(w_n)}+d^{(w_n)}$. From (\ref{b-speed}), (\ref{d-es}), and (\ref{I}),  there exists a constant  $\Tilde{C}_0>0$ depending only on $M_0,\delta_0$ and $f$ such that for every $x\in\R\backslash I^{\tau}_{\tau}$ and $t\in [0,\tau]$, it holds 
%\bel{bound-a}
%\big|a_{x}(t,x)\big|~\leq~\Tilde{C}_0\cdot (\tau-t)^{\alpha-1},\quad \big|a_{xx}(t,x)\big|~\leq~\Tilde{C}_0\cdot \left(1+|w_{xx}(t,x)|+{1\over |x|^{2-\alpha}}\right).
%\eeq
Multiplying the linear equation (\ref{v-lem42}) by $2v$, we have
\[
\big(v^2\big)_t + \big(a_n v^2\big)_x~=~ a_{n,x} v^2 ,\qquad v(0,x)~=~\bar{v}(x).
\]
Integrating the above equation over the domain
\[
\Omega~\doteq~ \bigcup_{t\in [0,\tau]} \{t\} \times (\R\backslash I^{\tau}_t)  ~=~\left\{(t,x)\in [0,\tau]\times\R: x\in \R\backslash I^{\tau}_t \right\}, 
\]
%
%\bel{I}
%I^{\tau}_t~\doteq~[-b_0(\tau-t)/2,b_0(\tau-t)/2]\,.
%\eeq
%
%, we get
%
%
%\textcolor{red}{
%{\bf Proof.}
%}
%\\ 
%1. Multiplying the linear equation (\ref{v-lem42}) by $2v$, we get
%\[
%\begin{cases}
%(v^2)_t + (a_n\, v^2)_x = a_{n,x}(t,x) \, v^2 \\
%v(0,x)~=~\bar{v}(x)
%\end{cases}
%\]
%Integrating the above equation over the domain 
%\[
%\Omega \doteq \Bigg\{(t,x) \,\,\Bigg| \,\, |x|>\frac{b_0}{2}(\tau - t), \, t\in [0,\tau] \Bigg\} 
%\]
%
%%\begin{figure}[htp]
%%        \centering
%%        \includegraphics[width=150mm,scale=0.3]{for-proof_4-2-2.jpg}
%%        \caption{The norm $\|v(\tau,\cdot)\|_{H^2(\R\backslash\{0\})}$ is estimated by using the balance laws for $v^2, v_x^2, v_{xx}^2$ on the domain $\Omega$.}
%%\end{figure}
and using the the first inequality in (\ref{bound-a}), we get 
\bel{v-est1}
\begin{split}
\int_{-\infty}^{\infty} v^2(\tau,x) dx&~\leq~ \int_{x\in \R\backslash I^{\tau}_0} \bar{v}^2(x) dx + \int_0^\tau\int_{x\in \R\backslash I^{\tau}_t} a_{n,x}(t,x) v^2(t,x)~dx dt\\
&\leq~\|\bar{v}(\cdot)\|^2_{L^2(\R\backslash I^{\tau}_0)} +\Tilde{C}_0 \int_0^{\tau} (\tau-t)^{\alpha-1}\cdot Z(t)~dt.
\end{split}
\eeq
%Thus, from (\ref{dx-a}), there exists a constant  $\Tilde{C}_0$ depending only on $M_0,\delta_0$ and $f$ such that 
%\[
%\|v(\tau,\cdot)\|^2_{L^2(\R\backslash\{0\})}~\leq~ \|\bar{v}(\cdot)\|^2_{L^2(\R\backslash I^{\tau}_0)} +\Tilde{C}_0\cdot \int_0^{\tau}\|v(t,\cdot)\|^2_{L^2(\R\backslash I^{\tau}_t)} \, dt.
%\]
%for some constant $\Tilde{C}_0$ depending only on $M_0,\delta_0$ and $f$. 
%By (\ref{dx-a}), taking the supremum of $|a_{n,x}(t,x)|$ over the set 
%\[
%\Omega_t \doteq \Bigg\{x \,\,\Bigg| \,\, |x|>\frac{b_0}{2}(\tau - t) \Bigg\}
%\]
%we obtain
%\[
%\|v(\tau,\cdot)\|^2_{L^2(\R\backslash\{0\})} \leq \|\bar{v}(\cdot)\|^2_{L^2(\Omega_0)} + \int_0^{\tau} \Big(b_1+\Gamma_f \left[C_0\ |x|^{\alpha-1}+4M_0\right]\Big) \, \|v(t,\cdot)\|^2_{L^2(\Omega_t)} \, dt
%\]
%
%
%
Similarly, differentiating equation (\ref{v-lem42}) with respect to $x$ and multiplying by $2v_x$ (and by $2v_{xx}$), we have
\[
\begin{cases}
\big(v_x^2\big)_t + \big(a_n\, v_x^2\big)_x ~=~ -a_{n,x} \, v^2_x,\qquad v_x(0,x)~=~\bar{v}'(x)\,,\\[4mm]
\big(v_{xx}^2\big)_t + \big(a_n\, v_{xx}^2\big)_x = -3a_{n,x} \, v_{xx}^2 - 2a_{n,xx} \, v_x \, v_{xx},\qquad v_{xx}(0,x)~=~\bar{v}''(x).
\end{cases}
\]
%\[
%\begin{cases}
%(v_x^2)_t + (a_n\, v_x^2)_x = -a_{n,x} \, v^2 \\
%v_x(0,x)~=~\bar{v}_x(x)
%\end{cases}
%\]
Integrating the above equations over $\Omega$, we  obtain 
\bel{v-est2}
\begin{split}
\int_{-\infty}^{\infty} v_x^2(\tau,x) dx&~\leq~\int_{x\in \R\backslash I^{\tau}_0} [\bar{v}']^2(x) \, dx- \int_0^\tau\int_{x\in \R\backslash I^{\tau}_t} a_{n,x}(t,x) v^2_x(t,x) ~dx dt\\
&~\leq~\|\bar{v}'(\cdot)\|^2_{L^2(\R\backslash I^{\tau}_0)} +\Tilde{C}_0 \int_0^{\tau} (\tau-t)^{\alpha-1} Z(t)~ dt,
\end{split}
\eeq
and
\bel{v-est3}
\begin{split}
\int_{-\infty}^{\infty} & v_{xx}^2(\tau,x) dx~\leq~\int_{x\in \R\backslash I^{\tau}_0} [\bar{v}'']^2(x) \, dx- \int_0^\tau\int_{x\in \R\backslash I^{\tau}_t} 3a_{n,x} v_{xx}^2 +2a_{n,xx} v_x v_{xx} ~dx dt\\
& \leq\|\bar{v}''(\cdot)\|^2_{L^2(\R\backslash I^{\tau}_0)} +3\Tilde{C}_0 \int_0^{\tau} (\tau-t)^{\alpha-1} Z(t) dt+2\Tilde{C}_0\int_0^\tau\int_{x\in \R\backslash I^{\tau}_t} a_{n,xx} v_x v_{xx} dx dt.
\end{split}
\eeq
Using the second inequality in (\ref{bound-a}), (\ref{I}), and H\"older's inequality, we estimate 
\begin{align*}
\int_{x\in \R\backslash I^{\tau}_t} a_{n,xx}& v_x v_{xx} dx~\leq~\int_{x\in \R\backslash I^{\tau}_t} \left(1+|w_{xx}(t,x)|+|x|^{\alpha-2}\right) \big|v_x(t,x) v_{xx}(t,x)\big| ~dx\\
&~\leq~\|v_x(t,\cdot)\|_{L^2(\R\backslash I^{\tau}_t)}\|v_{xx}(t,\cdot)\|_{L^2(\R\backslash I^{\tau}_t)}\\
& +\|v_x(t,\cdot)\|_{L^{\infty}(\R\backslash I^{\tau}_t)}\cdot \int_{x\in \R\backslash I^{\tau}_t}\big[|w_{n,xx}(t,x)|+|x|^{\alpha-2}\big]\cdot \big|v_{xx}(t,x)\big|~ dx \\
&~\leq~Z(t)+2Z^{1/2}(t)\cdot \int_{x\in \R\backslash I^{\tau}_t}\big[|w_{n,xx}(t,x)|+|x|^{\alpha-2}\big]\cdot \big|v_{xx}(t,x)\big|~dx \\
&~\leq~ Z(t)\cdot \left(1+2M_0+(2-\alpha)^{-1/2}\cdot  \Big[{b_0(\tau-t)\over 2}\Big]^{\alpha-3/2}\right).
\end{align*}
Thus, summing up (\ref{v-est1}), (\ref{v-est2}), and (\ref{v-est3}), we get
\[
Z(t) \leq Z(0)+\Tilde{C}_0 \int_{0}^{\tau}\left(2+4M_0+5(\tau-t)^{\alpha-1}+ \frac{2}{(2-\alpha)^{1/2}}  \Big[{b_0(\tau-t)\over 2}\Big]^{\alpha-3/2}\right) Z(t) dt,
\]
and by Gronwall's lemma, we derive for $\tau>0$ sufficiently small that 
\[
\|v(\tau,\cdot)\|^2_{H^2\bigl(\R\backslash\{0\}\bigr)}~=~Z(\tau)~\leq~ \frac{9}{4}\cdot Z(0)~=~\frac{9}{4}\cdot \|\bar{v}\|^2_{H^2
(\R\backslash I^{\tau}_0)},
\]
proving (\ref{32-ineq}).
\endproof
\medskip

As a consequence of Lemma \ref{Cr}, by Duhamel's formula, we get from (\ref{w-k+1}) that
\bel{est1-w-k+1}
\left\|w^{(k+1)}(\tau,\cdot)\right\|_{H^2(\R\backslash \{0\})}~\leq~{3\over 2}\cdot \|\overline{w}\|_{H^2
(\R\backslash I^{\tau}_0)}+{3\over 2}\cdot\int_{0}^{\tau}\left\|F^{(k)}(t,\cdot)\right\|_{H^2(\R\backslash I^{\tau}_t)}dt.
\eeq
for all $k\geq 1.$
\v

{\bf Step 2.} Using   Lemma \ref{F1}, Lemma \ref{Cr}, and (\ref{est1-w-k+1}),    we are now going to  establish  a priori estimates on  the sequence of approximations $\left(w^{(k)}\right)_{k\geq 1}$.
\begin{lemma}\label{Pr0} Assume that $w_n$ and $\vp_n$ satisfy (\ref{wn})-(\ref{nota-n}). There exists $T>0$ depending only on $M_0,\delta_0$, and $f$ such that{\color{blue}, for all $\tau \in [0,T]$,}
\bel{w-k-e1}
\left|w^{(k)}(\tau,0\pm)-\overline{w}(0\pm)\right|~\leq~{\delta_0\over 3},\qquad \left\|w^{(k)}(\tau,\cdot)\right\|_{H^2(\R\backslash\{0\})}~\leq~M_0.\eeq
In addition, the map $\tau\mapsto w^{(k)}(\tau,0\pm)$ is locally Lipchitz and  
\bel{w-k-e2}
\left|\dot{w}^{(k)}(\tau,0\pm)\right|~\leq~ 2\Gamma_1 \tau^{\alpha-1}~~~ a.e.~\tau\in [0,T].
\eeq
\end{lemma}
{\bf Proof.} {\bf 1.} Since $w^{(1)}(\tau,\cdot)=\overline{w}(\cdot)$ for all $\tau\in [0,T]$, the estimates (\ref{w-k-e1})-(\ref{w-k-e2}) hold for $k=1$. Assume that  (\ref{w-k-e1})-(\ref{w-k-e2}) hold for a given $k\geq 1$. Applying Lemma \ref{F1} for $w=w^{(k)}$, $\ell(t)=2\Gamma_1t^{\alpha-1}$ and $\delta=\ds {b_0(\tau-t)\over 2}$, we get 
\[
\begin{split}
\left\|F^{(k)}(t,\cdot)\right\|_{H^2(\R\backslash I^{\tau}_t)}&~\leq~\Gamma_1\cdot \left[t^{-3/4}+  (2\Gamma_1+1)t^{\alpha-1}\cdot \left({b_0(\tau-t)\over 2}\right)^{-3/4}\right],
\end{split}
\]
and (\ref{est1-w-k+1}), (\ref{cond}) yield
\bel{w-k-H}
\begin{split}
\left\|w^{(k+1)}(\tau,\cdot)\right\|_{H^2(\R\backslash \{0\})}&~\leq~{3M_0\over 4}+{3\over 2}\cdot \int_{0}^{\tau}\left\|F^{(k)}(t,\cdot)\right\|_{H^2(\R\backslash I^{\tau}_t)}dt\\
&~\leq~{3M_0\over 4}+\Gamma_{1,1}\cdot\left(\tau^{1/4}+\tau^{\alpha-3/4}\right)~\leq~{3M_0\over 4}+2\Gamma_{1,1}\cdot\tau^{\alpha-3/4}
\end{split}
\eeq
for some constant $\Gamma_{1,1}>0$ that depends on $\Gamma_1,b_0$ and $M_0$.
\medskip

{\bf 2.} For any $\tau\in [0,T]$ and $-\delta_1<\bar{x}_2<\bar{x}_1<0$,  consider the characteristics $t\mapsto x_i(t)\doteq x(t;\tau,\bar{x}_i)$ for $i=1,2$. Recalling Lemma \ref{F1}, Lemma \ref{Cr}, and (\ref{w-k+1}), we have 
\bel{w_x-k+1}
\begin{split}
\Big|&w^{(k+1)}(\tau,\bar{x}_2)-w^{(k+1)}(\tau,\bar{x}_1)\Big|\\
&\leq \left|\ov{w}\bigl(x_2(0)\bigr)-\ov{w}\bigl(x_1(0)\bigr)\right|+\int_{0}^{\tau}\left|F^{(k)}(t,x_2(t))-F^{(k)}(t,x_1(t))\right|dt\\
&\leq 2M_0K |\bar{x}_2-\bar{x}_1|+(2\Gamma^2_1+\Gamma_1)\int_{0}^{\tau}  \Big[ \Big({b_0t|\tau-t|\over 2}\Big)^{\alpha-1}+t^{\alpha-3/2}\Big] |x_2(t)-x_1(t)|dt\\
&\leq2M_0K |\bar{x}_2-\bar{x}_1|+\Gamma_{1,2}\cdot \left(\tau^{2\alpha-1}+\tau^{\alpha-1/2}\right)\cdot  |\bar{x}_2-\bar{x}_1|\\
&\leq\left(2M_0K+2\Gamma_{1,2}\cdot\tau^{1/4}\right)\cdot  |\bar{x}_2-\bar{x}_1|
\end{split}
\eeq
for some constant $\Gamma_{1,2}>0$ that depends on $\Gamma_1,b_0$. An entirely similar estimate holds for $\tau\in [0,T]$ and $0<\bar{x}_2<\bar{x}_1<\delta_1$.
\medskip

{\bf 3.} Given $0\leq \tau_1<\tau_2\leq T$, let $x^{\pm}_2(t)\doteq x(t;\tau_2,0\pm)$ be  the backward characteristic
starting from negative side and positive side of the  origin at time $\tau_2$ . From Lemma  \ref{F1}, Lemma \ref{Cr}, (\ref{w-k+1}), and (\ref{w_x-k+1}), one has
\[
\begin{split}
\Big|&w^{(k+1)}(\tau_2,0\pm)-w^{(k+1)}(\tau_1,0\pm)\Big|\\
&\leq~\Big|w^{(k+1)}(\tau_1,x^{\pm}_2(\tau_1))-w^{(k+1)}(\tau_1,0\pm)\Big|+\int_{\tau_1}^{\tau_2}\big|F^{(k)}(t,x^{\pm}_2(t))\big|dt\\
&\leq~6 M_0 K b_1 (\tau_2-\tau_1)+\Gamma_1 \int_{\tau_1}^{\tau_2}t^{\alpha-1}\left(1+ 2\Gamma_1\cdot\left[\left|{b_0(\tau-t)\over 2} \ln\left({b_0(\tau-t)\over 2}\right)\right|+t^{\alpha}\right] \right)dt\\
&\leq~6 M_0 K b_1 (\tau_2-\tau_1)+\left[\Gamma_1+2\Gamma^2_1\left(T^{\alpha}+\left|{b_0T\over 2} \ln\left({b_0T\over 2}\right)\right|\right)\right]\cdot\min\left\{{\tau_2^{\alpha}\over \alpha},\tau_1^{\alpha-1} (\tau_2-\tau_1)\right\}.
\end{split}
\]
In particular, for almost every $\tau\in [0,T]$, it holds  
\bel{dff}
\Big|w^{(k+1)}(\tau,0\pm)-\ov{w}(0\pm)\Big| ~\leq~ 6 M_0 K b_1 T + \left[\Gamma_1+2\Gamma^2_1\left(T^{\alpha}+\left|{b_0T\over 2} \ln\left({b_0T\over 2}\right)\right|\right)\right]\cdot {T^{\alpha}\over\alpha},
\eeq
and
\bel{w'-k+1} 
\Big|\dot{w}^{(k+1)}(\tau,0\pm)\Big| ~\leq~ 6 M_0 K b_1 + \left[\Gamma_1+2\Gamma^2_1\left(T^{\alpha}+\left|{b_0T\over 2} \ln\left({b_0T\over 2}\right)\right|\right)\right]\cdot \tau^{\alpha-1}.
\eeq
Thus, from (\ref{w-k-H}), (\ref{dff}), and (\ref{w'-k+1}), there exists a sufficiently small time $T>0$  depending only on $M_0,\delta_0$ and $f$ so that (\ref{w-k-e1})-(\ref{w-k-e2}) holds.
\endproof
\medskip

{\bf Step 3.} Thanks to the above estimates, we now complete step (i),  which is a key step toward the proof of Theorem \ref{t:main},  by  proving that the sequence of approximations
$w^{(k)}$ is Cauchy and converges to a solution $w$ of the linear problem (\ref{wne}).  
\begin{lemma}\label{w-n} Under the same settings  in Lemma \ref{Pr0}, the sequence of approximations $\big(w^{(k)}\big)_{k\geq 1}$ converges to a limit function $w$ in ${\bf L}^{\infty}([0,T],H^{2}(\R\backslash\{0\}))$ for  sufficiently small $T>0$  depending only on $M_0,\delta_0$, and $f$. Moreover, $w(\tau,\cdot)$ satisfies (\ref{w-k-e1})-(\ref{w-k-e2}) in $ [0,T]$.
\end{lemma}
{\bf Proof.} 
%To achieve the convergence of  $\big(w^{(k)}\big)_{k\geq 1}$ in ${\bf L}^{\infty}([0,T],H^{2}(\R\backslash\{0\}))$ for $T>0$ sufficiently small, we show that 
%\bel{sum-w}
%\sum_{k=1}^{\infty}\sup_{\tau\in [0,T]}\left\|w^{(k+1)}(\tau,\cdot)-w^{(k)}(\tau,\cdot)\right\|_{H^{2}(\R\backslash\{0\})}~<~\infty.
%\eeq
{\bf 1.} For every $k\geq 1$, we set  ${\bf w}^{(k)}\doteq w^{(k+1)}-w^{(k)}$ and  
\bel{wk-2-1}
\beta^{(k)}(\tau)~\doteq~ \sup_{t\in [0,\tau]}\left\|{\bf w}^{(k)}(t,\cdot)\right\|_{H^2(\R\backslash\{0\})}\qquad\forall \tau\in [0,T].
\eeq
By  (\ref{w-k-e2}), for almost every $\tau\in [0,T]$, we  define 
\bel{Sig-b}
\Sigma^{(k)}(\tau)~\doteq~ \max\left\{\left|\dot{\bf w}^{(k)}(\tau,0+)\right|,\left|\dot{\bf w}^{(k)}(\tau,0-)\right|\right\}~\leq~4\Gamma_1 \tau^{\alpha-1}.
\eeq
From (\ref{F2-1}) and (\ref{ln-wk}), ${\bf w}^{(k)}$ solves the semilinear equation
\[
v_t+a_n(t,x)\cdot v_x~=~{\bf F}^{\left(w^{(k)},w^{(k+1)}\right)}(t,x),\qquad v(0,\cdot)~=~0\,.
\]
In particular,  Lemma \ref{est-l1} and  the Duhamel formula yield
\[
\left\|{\bf w}^{(k+1)}(\tau,\cdot)\right\|_{H^2(\R\backslash\{0\})}~\leq ~ {3\over 2}\cdot \int_{0}^{\tau}\left\|{\bf F}^{\left(w^{(k)},w^{(k+1)}\right)}(s,\cdot)\right\|_{H^2(\R\backslash I^{\tau}_s)}ds
\]
Thus, by (\ref{I}), (\ref{w-k-e2}), and the second inequality in (\ref{F1-2}), we  derive 
\bel{beta-k}
\begin{split}
\beta^{(k+1)}(\tau)&\leq~C_1\cdot\Big(\beta^{(k)}(\tau) \int_{0}^{\tau}s^{\alpha-1} (\tau-s)^{\alpha-{3\over 2}} + (\tau-s)^{-{2\over 3}} s^{2\alpha-{11\over 6}}+\Sigma^{(k)}(s) (\tau-s)^{\alpha-{3\over 2}} \, ds\Big)\\
&\leq~C_2\beta^{(k)}(\tau)\tau^{2\alpha-{3\over 2}}+ C_1\cdot \int_{0}^{\tau}\Sigma^{(k)}(s) (\tau-s)^{\alpha-{3\over 2}} \, ds.
\end{split}
\eeq
for some constants $C_1,C_2$ depend only on $\Gamma_2$, and $b_0$
\medskip

{\bf 2.} Next, we are going to provide a  bound on  $\Sigma^{(k)}$ in terms of  $\beta^{(k)}$ and $\beta^{(k+1)}$. Given any $0<\tau_1<\tau_2\leq T$, let $\tilde{x}^{\pm}_{2}(t)=x_n(t;\tau_2,0\pm)$ be the characteristics, which reach the origin at time $\tau_i$, from the positive or negative side, respectively. From  (\ref{w-k+1}), it holds 
\[
{\bf w}^{(k+1)}(\tau_2,0\pm)~=~{\bf w}^{(k+1)}\left(\tau_1,\tilde{x}^{\pm}_2(\tau_1)\right)+\int_{\tau_1}^{\tau_2}{\bf F}^{\left(w^{(k)},w^{(k+1)}\right)}(t,\tilde{x}_2^\pm(t))dt.
\]
Using (\ref{x-n}), (\ref{wk-2-1}), (\ref{Sig-b}), and the first inequality in (\ref{F1-2}), we estimate 
\[
\begin{split}
\Big|{\bf w}^{(k+1)}(\tau_2,0\pm) & - {\bf w}^{(k+1)}(\tau_1,0\pm)\Big|\\
&\leq~\left|{\bf w}^{(k+1)}\left(\tau_1,\tilde{x}^{\pm}_2(\tau_1)\right)-{\bf w}^{(k+1)}(\tau_1,0\pm)\right|+\int_{\tau_1}^{\tau_2}\left|{\bf F}^{\left(w^{(k)},w^{(k+1)}\right)}\left(t,\tilde{x}_2^\pm(t)\right)\right|dt\\
& \leq~  \beta^{(k+1)}(\tau_1) \left|\tilde{x}_2^\pm(\tau_1)\right|+ 4\Gamma_2\Gamma_1\tau_1^{\alpha-1}\int_{\tau_1}^{\tau_2}  |\tilde{x}_2^\pm(t)|^\alpha \, dt \\
&\qquad\qquad\qquad\qquad~~~~~+ \Gamma_2  \beta^{(k)}(\tau_2)\int_{\tau_1}^{\tau_2} t^{\alpha-1} + 2\Gamma_1 t^{\alpha-1} \left(t^\alpha + |\tilde{x}_2^\pm(t)|^\alpha\right)dt\\
&\leq~C_3\cdot \left( \beta^{(k+1)}(\tau_1)+\tau_1^{\alpha-1} (\tau_2-\tau_1)^{\alpha}+\beta^{(k)}(\tau_2)\tau_1^{\alpha-1}\right)\cdot (\tau_2-\tau_1),\\
%& \le 2b_1 \beta^{(k+1)}(\tau_1) (\tau_2-\tau_1)+4\Gamma_2\Gamma_1 (2b_1)^\alpha\tau_1^{\alpha-1}\\
%& \quad + \Gamma_2 (2\Gamma_1+1) \beta^{(k)}(\tau_2) \Bigl(1 + \tau_2^\alpha + 2b_1 (\tau_2 - \tau_1)\Bigr) \int_{\tau_1}^{\tau_2} t^{\alpha-1} \, dt\\
%& \quad + \Gamma_2 \sup_{t \in [\tau_1,\tau_2]} \Sigma^{(k)}(t) \int_{\tau_1}^{\tau_2} [2b_1(\tau_2-\tau_1)]^\alpha \, dt\\
%& \le 2b_1 \beta^{(k+1)}(\tau_2) (\tau_2-\tau_1)\\
%& \quad + \Gamma_2 (2\Gamma_1+1) \beta^{(k)}(\tau_2) \tau_1^{\alpha-1} \Bigl(1 + \tau_2^\alpha + 2b_1 (\tau_2 - \tau_1)\Bigr) (\tau_2-\tau_1)\\
%& \quad + \Gamma_2 (2b_1)^\alpha \sup_{t \in [\tau_1,\tau_2]} \Sigma^{(k)}(t) (\tau_2-\tau_1)^{\alpha+1}
\end{split}
\]
and the increasing property of the map $\tau\mapsto  \beta^{(k+1)}(\tau) $ implies that
\begin{equation}\label{w-dot}
\left|\dot{\bf w}^{(k+1)}(s,0\pm)\right| ~\leq~ C_3 \cdot \Bigl(\beta^{(k+1)}(s) + \beta^{(k)}(s) s^{\alpha-1}\Bigr), \qquad ~a.e.~s \in ]0,T],
\end{equation}
where $C_3>0$ only depends on $\Gamma_1, \Gamma_2, b_1$, and $\alpha$.
\medskip

{\bf 3.} To complete the proof, we introduce the sequence of maps $\tau\mapsto \gamma^{(k)}(\tau)$ defined by 
\begin{equation*}
\gamma^k(\tau) ~ = ~ \frac{\beta^{(k)}(\tau)}{4 C_1}  + \int_0^\tau \Sigma^{(k)}(s) (\tau-s)^{\alpha-3/2} \, ds\qquad\forall \tau\in [0,T].
\end{equation*}
with $C_1$ being the same as in \eqref{beta-k}. Then, from \eqref{Sig-b}, \eqref{beta-k},  \eqref{w-dot}, and the increasing property of $\beta^{(k)}$, it holds
\[
\begin{split}
\gamma^{k+1}(\tau) &~ \le ~ \frac{C_2 \tau^{2\alpha-1/2}}{4 C_1} \beta^{(k)}(\tau) + \frac14 \int_0^\tau \Sigma^{(k)}(s) (\tau-s)^{\alpha-3/2} \, ds\\
&\qquad\qquad\qquad\qquad+C_3 \int_0^\tau \left(\beta^{(k+1)}(\tau) + \beta^{(k)}(\tau) s^{\alpha-1}\right) (\tau - s)^{\alpha - 3/2} \, ds\\
&~\leq~C_4 \tau^{\alpha-1/2} \beta^{(k+1)}(\tau) + C_4 \tau^{2\alpha - 3/2} \beta^{(k)}(\tau) + \frac14 \int_0^\tau \Sigma^{(k)}(s) (\tau-s)^{\alpha-3/2} \, ds
\end{split}
\]
for some constant $C_4>0$ that only depends on $C_1,C_2$ and $C_3$. In particular, by choosing $T>0$ sufficiently small such that 
\[
C_4 T^{\alpha-1/2}~\leq~{1\over 2},\qquad  C_4 T^{2\alpha - 3/2}~\leq~{1\over 16C_1},
\]
one obtains a contractive property %for  the sequence $\left(\gamma^{(k)}\right)_{k\geq 1}$ 
%\begin{multline*}
%\gamma^{k+1}(\tau) ~ \le ~ \frac{C_2 \tau^{2\alpha-1/2}}{4 C_1} \beta^{(k)}(\tau) + \frac14 \int_0^\tau \Sigma^{(k)}(s) (\tau-s)^{\alpha-3/2} \, ds\\
%+ C_3 \int_0^\tau \left(\beta^{(k+1)}(\tau) + \beta^{(k)}(\tau) s^{\alpha-1}\right) (\tau - s)^{\alpha - 3/2} \, ds\\
%\le ~ C_4 \tau^{\alpha-1/2} \beta^{(k+1)}(\tau) + C_4 \tau^{2\alpha - 3/2} \beta^{(k)}(\tau) + \frac14 \int_0^\tau \Sigma^{(k)}(s) (\tau-s)^{\alpha-3/2} \, ds,
%\end{multline*}
%for some constant $C_4>0$. Thus, for $\tau$ sufficiently small,
\begin{equation*}
\sup_{\tau\in [0,T]}\gamma^{(k+1)}(\tau) ~ \le ~ {1\over 2}\cdot  \sup_{\tau\in [0,T]} \gamma^{(k)}(\tau)\qquad\forall k\geq 1.
\end{equation*}
It follows that $\left(w^{(k)}\right)_{k\geq 1}$ is  a Cauchy sequence in ${\bf L}^{\infty}([0,T],H^{2}(\R\backslash\{0\}))$ and converges to a function $w_{n+1} \in {\bf L}^{\infty}([0,T],H^{2}(\R\backslash\{0\}))$. This  implies that 
\[
\|w_{n+1}(\tau,\cdot)\|_{H^2(\R\backslash \{0\})}~\leq~M_0,\qquad \lim_{k\to\infty}w^{(k)}(\tau,0\pm)~=~w_{n+1}(\tau,0\pm)\qquad\forall \tau\in [0,T],
\]
and $w_{n+1}$ satisfies (\ref{w-k-e1})--(\ref{w-k-e2}).
Furthermore, from (\ref{w-dot}) and (\ref{w-k-e2}), the limit  $\ds\lim_{k\to\infty}\dot{w}^{(k)}(\tau,0\pm)$ exist and bounded by $2\Gamma_1 \tau^{\alpha-1}$ for almost every $\tau\in (0,T)$. Thus, for every $0<\tau_1<\tau_2<T$, one has 
\[
\begin{split}
w_{n+1}(\tau_2,0\pm)-w_{n+1}(\tau_1,0\pm)&=~\lim_{k\to\infty} \left(w^{(k)}(\tau_2,0\pm)-w^{(k)}(\tau_1,0\pm)\right)\\
&=~\lim_{k\to\infty}\int_{\tau_1}^{\tau_2}\dot{w}^{(k)}(\tau,0\pm)d\tau~=~\int_{\tau_1}^{\tau_2}\lim_{k\to\infty}\dot{w}^{(k)}(\tau,0\pm)d\tau\,,
\end{split}
\]
and this yields
\[
\lim_{k\to\infty}\dot{w}^{(k)}(\tau,0\pm)~=~\dot{w}_{n+1}(\tau, 0\pm)\quad a.e.~\tau\in [0,T].
\]
Thus, recalling  the first estimate in (\ref{F1-2}), we have 
\[
\lim_{k\to\infty}F\big(\tau,x,w^{(k)}(\tau,x)\big)~=~F\big(\tau,x,w_{n+1}(\tau,x)\big).%\quad\quad a.e.~\tau\in [0,T], x\neq 0.
\]
Finally, taking $k\to+\infty$ in (\ref{w-k+1}), we obtain that for all $t_0\in [0,T]$,
\begin{equation*}
w_{n+1}(t_0,x_0) ~ = ~ \overline{w}(x_n(0;t_0,x_0)) + \int_0^{t_0} F\Bigl(t,x_n(t;t_0,x_0),w_{n+1}\Bigr) \, dt,
\end{equation*}
and $w_{n+1}$ is a solution to the semilinear equation (\ref{BH-3}).
\endproof
\subsection{Proof of Theorem \ref{t:main}}
{\bf 1.} By Lemma \ref{w-n}, we  inductively construct the sequence of approximate solutions $(w_n)_{n\geq 1}$ where each $w_n$ solves (\ref{wne})  and satisfies (\ref{wn}), on a suitably small time interval $[0,T]$. Moreover, the map $\tau\mapsto w_n(\tau,0\pm)$ is locally Lipschitz and  
\[
\left|\dot{w}_n(\tau,0\pm)\right|~\leq~ 2\Gamma_1\cdot \tau^{\alpha-1}~~~ a.e.~\tau\in [0,T].
\]
As outlined at the end of Section \ref{sec:main}, we show that the sequence $(w_n)_{n\geq 1}$  is Cauchy w.r.t. the norm in $H^1(\R\backslash\{0\})$, hence it converges to a unique limit $w$ providing an entropic solution to the Cauchy problem (\ref{BH-2}). In order to do so, for a fixed $n\geq 2$, we define 
\[
\begin{cases}
{\bf w}_n~\doteq~w_n-w_{n-1},\qquad { \Psi}_n~\doteq~\vp^{(w_n)}-\vp^{(w_{n-1})}, \qquad {\bf u}_n~\doteq~{\bf w}_n+{ \Psi}_n,\\[4mm]
A_n~\doteq~a_{n}-a_{n-1},\qquad \beta_n(\tau)~\doteq~\ds\sup_{t\in [0,\tau]}\left\|{\bf w}_{n}(t,\cdot)\right\|_{H^1(\R\backslash\{0\})}.
\end{cases}
\]
Recalling {\color{blue} (\ref{usw}) and} (\ref{wne}), we have 
\[
{\bf u}_{n+1,t}+ a_n\cdot {\bf u}_{n+1,x}~=~{\bf G}[{ \Psi}_{n+1}+{\bf w}_{n+1}]-A_nw_{n,x}-A_{n+1}\vp^{(w_{n+1})}_x.
\]
From Lemma \ref{est-l1} and Duhamel's formula, we obtain that for all $\tau\in [0,T]$, we have
\begin{align}
\|{\bf u}_{n+1}(&\tau,\cdot)\|_{H^1(\R\setminus\{0\})} \nonumber \\
& \leq ~ \frac32 \int_0^\tau \left\|\Bigl({\bf G}[\Psi_{n+1} + {\bf w}_{n+1}] - A_n w_{n,x} - A_{n+1} \vp_x^{(w_{n+1})}\Bigr)(t,\cdot)\right\|_{H^1(\R\setminus I^\tau_t)} dt \label{un+1}
\end{align}

\medskip

{\bf 2.} To bound the right-hand side of (\ref{un+1}), we first recall the last inequality in (\ref{Psi-2-1}) and  Lemma \ref{B-est} to get 
\[
\begin{cases}
\|{\bf G}[\Psi_{n+1}\|_{H^1(\R\setminus \{0\})}~\leq~\Gamma_{f,1}\cdot \|{\bf w}_{n+1}(t,\cdot)\|_{H^1(\R\setminus\{0\})} ~ \le ~ \Gamma_{f,1}\cdot \beta_{n+1}(t),\\[4mm]
\|{\bf G}[{\bf w}_{n+1}\|_{H^1(\R\setminus[-\delta,\delta])}~\leq~\ds{\|{\bf w}_{n+1}(t,\cdot)\|_{H^1(\R\setminus\{0\})}\over \delta^{1/2}}~\leq~{\beta_{n+1}(t)\over \delta^{1/2}},
\end{cases}
\]
for $\delta>0$ sufficiently small. On the other hand, since  $a_n=b^{(w_n)}+d^{(w_n)}$, from (\ref{b-speed-2-1}) and (\ref{d-w}), it holds
\bel{An}
\begin{cases}
|A_n(t,x)| \leq \Tilde{C}_1 \|{\bf w}_{n}(t,\cdot)\|_{H^1(\R\setminus\{0\})},\\[4mm]
 |A_{n,x}(t,x)|\leq \Tilde{C}_1 \left({\|{\bf w}_{n}(t,\cdot)\|_{H^1(\R\setminus\{0\})}\over |x|^{1-\alpha}} + {\bf w}_{n,x}(t,x)\right)
 \end{cases}
\eeq
and this implies that
\[
\|A_n(t,\cdot)\|_{H^1(\R\backslash \{0\})} ~ \leq ~\Tilde{C}_2\cdot \|{\bf w}_{n}(t,\cdot)\|_{H^1(\R\setminus\{0\})},
\]
for some constant $\Tilde{C}_1,\Tilde{C}_2>0$  depending on $f$, $M_0$, and $C_0$. Thus, using (\ref{w-1-b}) and (\ref{vpb1})-(\ref{vpb2}), we derive
%\[
%\begin{split}
%\|A_n(t,\cdot) w_{n,x}(t,\cdot)\|_{H^1(\R\setminus \{0\})}~\leq~ \|A_n(t,\cdot) \|_{H^1(\R\setminus \{0\})}\cdot \|w_n\|_{H^1(\R\setminus\{0\})}
%\end{split}
%\]
\begin{align*}
\|A_n(t,\cdot) w_{n,x}(t,\cdot)\|_{H^1(\R\setminus \{0\})}&\leq 4 \|A_n(t,\cdot) \|_{H^1(\R\setminus \{0\})}\cdot \|w_{n,x}\|_{H^1(\R\setminus\{0\})} \\[4mm]
&\le 4M_0\cdot \|A_n(t,\cdot) \|_{H^1(\R\setminus \{0\})} \leq 4\Tilde{C}_2M_0 \beta_n(t),\\[4mm]
\ds \left\|A_{n+1}(t,\cdot) \vp_x^{(w_{n+1})}(t,\cdot)\right\|_{H^1(\R\setminus[-\delta,\delta])} &\le 4\|A_{n+1}(t,\cdot) \|_{H^1(\R\setminus \{0\})}\|\vp_x^{(w_{n+1})}(t,\cdot)\|_{H^1(\R\setminus[-\delta,\delta])} \\[4mm]
&\le \ds \frac{4C_0}{\delta^{3/2-\alpha}}\|A_{n+1}(t,\cdot) \|_{H^1(\R\setminus \{0\})} \le \frac{4C_0\Tilde{C}_2\beta_{n+1}(t)}{\delta^{3/2-\alpha}} 
\end{align*}
Thus, from (\ref{un+1}) and (\ref{I}), we obtain
\begin{align}\label{u-n+1-e}
& \|{\bf u}_{n+1}(\tau,\cdot)\|_{H^1(\R\setminus\{0\})} \\
& \leq {3\over 2}\int_{0}^{\tau} 4\Tilde{C}_2M_0\beta_n(t)+\left(\Gamma_{f,1}+\left|{2\over b_0(\tau-t)}\right|^{1/2}+4C_0\Tilde{C}_2\left|{2\over b_0(\tau-t)}\right|^{3/2-\alpha}\right)\beta_{n+1}(t)dt \nonumber\\
%&~\le ~ \frac32 \int_0^\tau \left\|\Bigl({\bf G}[\Psi_{n+1} + {\bf w}_{n+1}] - A_n w_{n,x} - A_{n+1} \vp_x^{(w_{n+1})}\Bigr)(t,\cdot)\right\|_{H^1(\R\setminus I^\tau_t)} \, dt\\
&\le 6\Tilde{C}_2M_0\tau \beta_{n}(\tau)+\Tilde{C}_3 \tau^{\alpha-1/2}\beta_{n+1}(\tau) \nonumber
\end{align}
for some constant $\Tilde{C}_3>0$ depending on $C_0,\Tilde{C}_2,\Gamma_{f,1},b_0$ and $\alpha$. From the last inequality in (\ref{Psi-2-1}), one has 
\[
\begin{split}
\|{\bf u}_{n+1}(\tau,\cdot)\|_{H^1(\R\setminus\{0\})}&~\geq~\|{\bf w}_{n+1}(\tau,\cdot)\|_{H^1(\R\setminus\{0\})}-\|\Psi_{n+1}(\tau,\cdot)\|_{H^1(\R\setminus\{0\})}\\ 
&~ \geq ~ \left(1-\Gamma_{f,1}  \tau^{\alpha-1/2}\right)\cdot\|{\bf w}_{n+1}(\tau,\cdot)\|_{H^1(\R\setminus\{0\})},
\end{split}
\]
and (\ref{u-n+1-e})  yields 
\[
\begin{split}
\|{\bf w}_{n+1}(\tau,\cdot)\|_{H^1(\R\setminus\{0\})} &~ \le ~ {\tau^{\alpha-1/2}\over 1-\Gamma_{f,1}  \tau^{\alpha-1/2}}\cdot \left({3\over 2}\Tilde{C}_2M_0\tau^{3/2-\alpha}\cdot \beta_{n}(\tau)+\Tilde{C}_3\cdot\beta_{n+1}(\tau)\right)\\
&~\leq~{1\over 4}\cdot \beta_n(\tau)+{1\over 2}\cdot \beta_{n+1}(\tau)
\end{split}
\]
for all  $\tau\in [0,T]$ with $T>0$ sufficiently small. In particular,
\[
\beta_{n+1}(\tau) ~\leq~\frac12 \cdot \beta_n(\tau) \qquad \text{for all } \tau \in [0,T],
\]
 and $\left(w_{n}\right)_{n\geq 1}$ is a Cauchy sequence in ${\bf L}^{\infty}([0,T],H^{1}(\R\backslash\{0\}))$ which  converges to the unique limit $w$ such that  
 \[
w(t,0-)-w(t,0+)~\geq~{\delta_0\over 3},\qquad \|w(t,\cdot)\|_{H^2(\R\backslash\{0\})}~\leq~M_0,\quad\forall t\in [0,T].
\]
Moreover, the map $t\mapsto w(t,0\pm)$ is locally Lipschitz and 
\[
\big|\dot{w}(t,0\pm)|~\leq~2\Gamma_1\cdot t^{\alpha-1}\quad a.e.~ t\in [0,T].
\]
Hence,  $u\doteq w+\vp^{(w)}$ satisfies (i)-(ii) in Definition \ref{d:12} and (\ref{Lp-u}). To verify that $u$ a piecewise regular solution to (\ref{BL-1})-(\ref{IC1}), we notice that $w_{n+1}$ is the solution to (\ref{wne}) and 
\[
{d\over dt}\Big[w_{n+1}+\vp^{(w_{n+1})}\Big]+a_n\cdot {d\over dx}\Big[w_{n+1}+\vp^{(w_{n+1})}\Big]~=~{\bf G}\left[w_{n+1}+\vp^{(w_{n+1})}\right]-A_{n+1}\cdot \vp^{(w_{n+1})}_x.
\]
Denoting by $x_n(\cdot;t_0,x_0)$ the solution to
\[
\dot{x}(t) ~ = ~ a_n(t,x(t)), \qquad x(t_0) ~ = ~ x_0,
\]
the formula above implies
\bel{wn+1-c}
\begin{split}
\Big[& w_{n+1}+\vp^{(w_{n+1})}\Big](t_0,x_0)~=~(\ov{w}+\ov{v})(x_n(0;t_0,x_0))+\\
& +\int_{0}^{t_0}{\bf G}\left[w_{n+1}+\vp^{(w_{n+1})}\right](t,x_n(t;t_0,x_0))dt -\int_{0}^{t}\big[A_{n+1} \vp^{(w_{n+1})}_x\big](t,x_n(t;t_0,x_0))dt
\end{split}
\eeq
From the first inequality in (\ref{An}), it holds
\begin{equation*} 
\begin{split}
\lim_{n\to\infty} \int_{0}^{t}\big[& A_{n+1}\cdot \vp^{(w_{n+1})}_x\big](t,x(t;t_0,x_0))~dt\\
& \leq ~ {\color{blue}\lim_{n\to\infty}} \Tilde{C}_1\cdot \beta_{n+1}(t_0)\int_{0}^{t_0}\left|\vp^{(w_{n+1})}(t,x(t;t_0,x_0))\right|~dt~=~0.
\end{split}
\end{equation*}
Taking $n\to\infty$ in (\ref{wn+1-c}), we obtain 
\[
u(t_0,x_0)~=~(\ov{w}+\ov{v})(x(0;t_0,x_0))+\int_{0}^{t_0}{\bf G}\left[u(t,\cdot)\right](x(t;t_0,x_0))~dt
\]
with $t\mapsto x(t;t_0,x_0)$ being  the characteristics curve, obtained by solving
\[
\dot{x}~=~\tilde{a}(t,x,u)~\doteq~\left(f'(u(t,x))- {f(u^-(t))-f(u^+(t))\over u^{-}(t)-u^{+}(t)}\right),\qquad x(t_0)~=~x_0.
\]

{\bf 3.} It remains to prove the uniqueness of (\ref{BL-1})-(\ref{IC1}). Assume that $\tilde{u}$ is a piecewise regular  solution to  (\ref{BL-1})-(\ref{IC1}). Then we have
\[
\sup_{t\in [0,T]}\|\tilde{u}(t,\cdot)\|_{H^1(\R\backslash\{0\})}~\doteq~M_1~<~\infty,
\]
and for every $\delta>0$ there exists a constant $M_{\delta}>0$ such that 
\[
|\tilde{u}_x(t,x)|~\leq~M_{\delta}\qquad\forall (t,x)\in [0,T]\times \R\backslash (-\delta,\delta).
\]
By  (\ref{RH1}), (\ref{cc2}), and the continuity ${\bf G}[\tilde{u}(t,\cdot)](\cdot)$ outside the origin, $\tilde{u}$ is continuously differentiable with respect to both variables $t,x$ for $x\neq 0$. In particular, the map $t\mapsto \tilde{u}_x(t,0\pm)$ is continuous and 
\[
\inf_{t\in [0,T]}\tilde{u}(t,0-)-\tilde{u}(t,0+)~>~0.
\]
Following the same argument as in the proof of Lemma \ref{Cr}, there exist $\tilde{\delta}_1,\tilde{b}_0>0$ and $0<T_1\leq T$ small such that for all $\tau\in [0,T_1]$ and $t\in [0,\tau]$, it holds 
\[
\tilde{a}(t,x,\tilde{u})\cdot\mathrm{sign}(x)~<~0\qquad\forall x\in [-\tilde{\delta}_1,\tilde{\delta}_1]\backslash\{0\},
\]
and
\bel{x-t-p}
\begin{cases}
\tilde{x}(t;\tau,x_{\tau})~\geq~x_{\tau}+\tilde{b}_0\cdot (\tau-t)&\mathrm{if}~~x_{\tau}\in ]0,\tilde{\delta}_1]\,,\\[4mm]
\tilde{x}(t;\tau,x_{\tau})~\leq~x_{\tau}-\tilde{b}_0\cdot (\tau-t)&\mathrm{if}~~ x_{\tau}\in [-\tilde{\delta}_1,0[\,,
\end{cases}
\eeq
where $\tilde{x}(\cdot;\tau,x_\tau)$ is the solution to
\[
\dot{x}(t) ~ = ~ \tilde{a}(t,x(t),\tilde{u}), \qquad \tilde{x}(\tau) ~ = ~ x_\tau.
\]
%{\bf 5.} 
%
%and this implies that 
%\bel{sign-ta}
%\tilde{a}(t,x,\tilde{u}(t,x))\cdot \sign(x) ~<~0\quad\forall (t,x)\in [0,T_1]\times \big([-\delta_1,\delta_1]\backslash \{0\}\big)
%\eeq
%for some constant $\delta_1>0$.
%\medskip
%$\tilde{\delta}_1,\tilde{b}_0, \tilde{b}_1>0$ and $0<T_2<T_1$ small such that for all $t\in [0,T_1]$
%\bel{xi-i}
%\tilde{\xi}(t;\tau,x_\tau)~\notin~\tilde{I}^{\tau}_t,\qquad \tilde{b}_0(\tau-t)~\leq~\big|\tilde{\xi}(t;\tau,x_\tau)-x_\tau\big|~\leq~\tilde{b}_1(\tau-t)
%\eeq
%with 
%\bel{I}
%\tilde{I}^{\tau}_t~\doteq~[-\tilde{b}_0(\tau-t),\tilde{b}_0(\tau-t)]\,.
%\eeq
%
For every $0<\delta<\tilde{\delta}_1$ small, we shall provide a better upper bound for  
\[
\tilde{\gamma}_{\tau}(\delta,t)~=~\sup_{x\in \R\backslash [\tilde{x}(t;\tau,-\delta),\tilde{x}(t;\tau,\delta)]}~|u_x(t,x)|\qquad\forall t\in [0,\tau].
\]
For every  $x_1<x_2<\tilde{x}(t;\tau,-\delta)$ or $\tilde{x}(t;\tau,\delta)<x_1<x_2$, set $\tilde{x}_i(\cdot)\doteq \tilde{x}(\cdot;t,x_i)$ for $i=1,2$. Then, one has
\[
\begin{split}
\left|\dot{\tilde{x}}_1(s)-\dot{\tilde{x}}_2(s)\right|&~\leq~ \left|\tilde{a}(s,x_1(s),\tilde{u})-\tilde{a}(s,x_2(s),\tilde{u})\right|\\
&~\leq~\Tilde{C_2}\cdot |\tilde{u}(s,x_1(s))-\tilde{u}(s,x_2(s))|~\leq~\Tilde{C}_2\cdot \tilde{\gamma}_{\tau}(\delta,s)\cdot \left|\tilde{x}_1(s)-\tilde{x}_2(s)\right|
\end{split}
\]
with $\Tilde{C}_2=\max_{\omega\in [-2M_1,2M_1]}\big|f''(\omega)\big|$. Applying Gronwall's inequality, we get 
\[
\left|\tilde{x}_1(s)-\tilde{x}_2(s)\right|~\leq~\exp \left(\Tilde{C}_2\cdot \int_{s}^{t}\tilde{\gamma}_{\tau}(\delta,r)dr\right)\cdot |x_1-x_2|\quad\forall s\in [0,t].
\]
%and this yields 
%\[
%|\tilde{\xi}(t;\tau,x_1)-\tilde{\xi}(t;\tau,x_2)|~\leq~\left(\exp{\int_{t}^{\tau}}\O(1)\cdot\tilde{\gamma}_{\delta}(s)ds\right)\cdot |x_2-x_1|\qquad\forall t\in [0,\tau].
%\]
%Notice that 
%\[
%\tilde{\xi}(t;\tau,x_1), \tilde{\xi}(t;\tau,x_2)~\leq~-\delta-\tilde{b}_0(\tau-t) \qquad\forall t\in [0,\tau].
%\]
By Lemma \ref{B-est}, (\ref{Xadef}), and (\ref{x-t-p}), we estimate 
\[
\begin{split}
|\tilde{u}(t,x_1)&-\tilde{u}(t,x_2)|~\leq~\big|(\bar{v}+\ov{w})(\tilde{x}_1(0))-(\bar{v}+\ov{w})(\tilde{x}_2(0))\big|\\
&\qquad\qquad\qquad\qquad+\int_{0}^{t}\left|{\bf G}[\tilde{u}(t,\cdot)](\tilde{x}_1(s))-{\bf G}[\tilde{u}(s,\cdot)](\tilde{x}_2(s))\right|ds\\
&\leq~\big(2M_0+\Tilde{C}_3\cdot (\delta+\tau)^{\alpha-1}\big)\cdot |\tilde{x}_1(0)-\tilde{x}_2(0)|\\
&\qquad\quad+\Tilde{C}_3\cdot M_1\int_{0}^{t}\left(\ln^2\big[\delta+\tilde{b}_0(\tau-s)\big]+{1\over \delta+\tilde{b}_0(\tau-s)}\right)\cdot |\tilde{x}_2(s)-\tilde{x}_1(s)| ds\\
&\leq~\Tilde{C}_4\cdot \big((\delta+\tau)^{\alpha-1}+\big|\ln (\delta+\tilde{b}_0(\tau-t)) \big|\big)\cdot \exp \left(\Tilde{C}_2\cdot \int_{0}^{t}\tilde{\gamma}_{\tau}(\delta,s)ds\right)\cdot |x_1-x_2|.
\end{split}
\]
Thus, for $t\in [0,\tau]$, we have 
\bel{ODE-gam}
\tilde{\gamma}_{\tau}(\delta,t)~\leq~\Tilde{C}_4\cdot \big((\delta+\tau)^{\alpha-1}+\big|\ln (\delta+\tilde{b}_0(\tau-t)) \big|\big)\cdot \exp \left(\Tilde{C}_2\cdot \int_{0}^{t}\tilde{\gamma}_{\tau}(\delta,s)ds\right)
\eeq
Equivalently,  
\[
-{d\over dt}\exp \left(-\Tilde{C}_2\cdot \int_{0}^{t}\tilde{\gamma}_{\tau}(\delta,s)ds\right)~\leq~\Tilde{C}_2\Tilde{C}_4\cdot \big((\delta+\tau)^{\alpha-1}+\big|\ln (\delta+\tilde{b}_0(\tau-t)) \big|\big),
\]
and thus there exists a small time $0<T_2<T_1$ and a constant $\Tilde{C}_5>0$ such that
\[
\exp \left(\Tilde{C}_2\cdot \int_{0}^{\tau}\tilde{\gamma}_{\tau}(\delta,s)ds\right)~\leq~\Tilde{C}_5\qquad\forall \tau\in [0,T_2].
\]
Recalling (\ref{ODE-gam}), we finally get 
\bel{bb-1}
\sup_{x\in \R\backslash [-\delta,\delta]}~|\tilde{u}_x(\tau,x)|~=~\tilde{\gamma}_{\tau}(\delta,\tau)~\leq~\Tilde{C}_6\cdot \delta^{\alpha-1}\quad\forall \tau\in [0,T_2]
\eeq
for some constant $\Tilde{C}_6>0$ which does not depend on $\delta$.
\medskip
%
%\textcolor{blue}{\bf TO DO:} \textcolor{blue}{\bf Please check everything above.}
%\medskip

{\bf 4.} Finally, to show that $\tilde{u}(t,\cdot)$ coincides with $u(t,\cdot)$ for all  $t\in [0,T]$,   defining
$$
{\bf u}(t,x)\doteq \tilde{u}(t,x)-u(t,x),\qquad A(t,x) ~ \doteq~ \tilde{a}(t,x,\tilde{u}) - \tilde{a}(t,x,u),
$$
we have  
\begin{equation*}
{\bf u}_t + \tilde{a}(t,x,\tilde{u})\cdot {\bf u}_x ~=~ {\bf G}[{\bf u}] - A(t,x)\cdot u_x.
\end{equation*}
Multiplying the above  equation  by $2{\bf u}$, we derive 
\bel{bf-u}
\big({\bf u}^2\big)_t+\big(\tilde{a}(t,x,\tilde{u})\cdot {\bf u}^2\big)_x~=~{\bf u}\cdot \big[2{\bf G}[{\bf u}]-2A(t,x)u_x+ \tilde{a}_x(t,x,\tilde{u}) {\bf u}\big].
\eeq
For every $\tau\in [0,T_2]$, integrating (\ref{bf-u}) over the domain 
$$\Omega_{\tau}~\doteq~\big\{(t,x)\in [0,\tau]\times \R: x\in \R\backslash [\tilde{x}_{\tau}^{-}(t),\tilde{x}_{\tau}^{+}(t) ]\big\},$$
where $\tilde{x}_\tau^\pm(t) \doteq \tilde{x}(t;\tau,0\pm)$, we get  
\[
\begin{split}
\|{\bf u}(\tau,\cdot)\|^2_{L^2(\R)}&~\le~ \int_{0}^{\tau} \int_{\R\backslash [\tilde{x}_{\tau}^{-}(t),\tilde{x}_{\tau}^{+}(t)]} \left|{\bf u}\cdot \big[2{\bf G}[{\bf u}]-2A(t,x)u_x+ \tilde{a}_x(t,x,\tilde{u}) {\bf u}\big]\right| dxdt.
\end{split}
\]
Similarly, arguing as in the proof of Lemma \ref{est-l1}, we obtain
\begin{multline*}
\|{\bf u}_x^2(\tau,x)\|^2_{L^2(\R)}\\
\leq~\int_0^\tau \int_{\R\backslash [\tilde{x}_{\tau}^{-}(t),\tilde{x}_{\tau}^{+}(t)]} {\bf u}_x \left[2\frac{d}{dx} {\bf G}[{\bf u}] - \tilde{a}_x(t,x,\widetilde{u}) {\bf u}_x - 2A_x(t,x) u_x - 2A(t,x) u_{xx}\right] dx dt.
\end{multline*}
From Lemma \ref{B-est}, (\ref{bb-1}), the two inequalities above, and the fact that $\alpha > 3/4$, there exists $T_3>0$ so small that for every $\tau \in [0,T_3]$
\begin{equation*}
\|{\bf u}(\tau,\cdot)\|_{H^1(\R \setminus \{0\})}^2 ~ \le ~ \widetilde{C}_7 \int_0^\tau [\chi_\tau(t)]^{-1/2} \|{\bf u}(t,\cdot)\|_{H^1(\R \setminus \{0\})}^2 dt,
\end{equation*}
where $\widetilde{C}_7 > 0$ does not depend on $\tau$ and $\chi_\tau \doteq \max\{\tilde{x}_\tau^+,-\tilde{x}_\tau^-\}$. From (\ref{x-t-p}) and Gronwall's inequality, we conclude that
\[
\|{\bf u}(\tau,\cdot)\|_{H^1(\R\backslash\{0\})}~=~0\qquad\forall \tau\in [0,T_3].
\]
Finally, we set
\[
\Tilde{T}~\doteq~\sup~\{\tau\in [0,T]:\tilde{u}(t,\cdot)=u(t,\cdot)\qquad\forall t\in [0,\tau]\}.
\]
By the continuity of $\tilde{u},u$ outside the origin,   $\tilde{u}(\widetilde{T},\cdot) = u(\widetilde{T},\cdot)$ has the same regularity as $\overline{w} + \overline{v}$. Consequently, if $\widetilde{T} < T$, then arguing as above we can find $\overline{T} \in \, ]\widetilde{T},T]$ such that $\tilde{u}(\tau,\cdot) = u(\tau,\cdot)$ for every $\tau \in [0,\overline{T}]$, which contradicts the definition of $\widetilde{T}$.
\endproof
\appendix
\section{Estimates on the nonlocal source} 
\setcounter{equation}{0}
In this section, we shall establish basic estimates which are used in the proof of Lemma \ref{F1} and Lemma \ref{F2}. Assume that $K$  satisfies  {\bf (H1)}-{\bf (H2)}. Recalling the definition of $\Lambda$, $\Phi$, $\eta$ and $\phi$ in (\ref{Lambda_1})-(\ref{Lambda_2}) and (\ref{Phi})-(\ref{phi-1}), we set 
\[
\phi_b(x)~\doteq~\phi(x,b)~=~\eta(x)\cdot [\Phi(b)-\Phi(x+b)]~=~\eta(x)\cdot\left[\int_{0}^{b}\Lambda(y)\,dy-\int_{0}^{x+b}\Lambda(y)\,dy\right].
\]
We first provide some bounds on ${\bf G}[\chi_{[0,\infty[}\phi_b]$ with $\chi_{[0,\infty[}$ being the indicator function on $[0,\infty[$ for $b>0$ small. As usual, by the Landau symbol $\O(1)$ we shall denote a uniformly bounded quantity which does not depend on $b$.
\begin{lemma}\label{A1} Assume that $0<b<1/4$.  For every $0<|x|<1/4$ and $\delta>0$, we have 
\bel{Gp-1}
\begin{cases}
\left|{\bf G}[\chi_{[0,\infty[}\phi_b](x)\right|~\leq~\O(1),\qquad \left|\ds{d\over dx}{\bf G}[\chi_{[0,\infty[}\phi_b](x)\right|~\leq~\O(1)\cdot\ln^2|x|,\\[4mm]
\left|\ds{d^2\over dx^2}{\bf G}[\chi_{[0,\infty[}\phi_b](x)\right|~\leq~\ds\O(1)\cdot \left|{\ln |x|\over x}\right|,\qquad \left\|{\bf G}[\chi_{[0,\infty[}\phi_b]\right\|_{H^2(\R\backslash [-\delta,\delta])}~\leq~\O(1)\cdot \delta^{-2/3}.
\end{cases}
\eeq
\end{lemma}
{\bf Proof.}  From Lemma \ref{G1} and $\phi_b(0)=0$, it holds
\[
{\bf G}\big[\chi_{[0,\infty[}\phi_b\big](x)~=~-\int_{0}^{2}\Phi'(y+b)\cdot\Lambda(x-y)~dy\qquad\forall |x|< 1/2.
\]
\noindent {\bf Case 1.} Assume that $-1/4<x<0$. Observe that for all $y\in [-2,2]\backslash\{0\}$,
\[
|\Lambda(y)|, |\Phi'(y)|~\leq~\O(1)\cdot \left(\big|\ln|y|\big|+1\right),\qquad |\Phi(y)|~\leq~\O(1)\cdot |y|\cdot\left(\big|\ln|y|\big|+1\right),
\]
we have 
\begin{multline*}
\left|{\bf G}\big[\chi_{[0,\infty[}\phi_b\big](x)\right|~=~\left|\int_{0}^{2}\Phi'(y+b)\cdot\Lambda(x-y)~dy\right|\\~\leq~\O(1)\cdot\int_{0}^{2}\big(1+\big|\ln|y+b|\big)\big(1+\big|\ln |x-y|\big|\big)~dy~\leq~\O(1).
\end{multline*}
%\[
%\left|{\bf G}\big[\chi_{[0,\infty[}\phi_b\big](x)\right|~=~\left|\int_{0}^{2}\Phi'(y+b)\cdot\Lambda(x-y)~dy\right|~\leq~\O(1)\cdot\int_{0}^{2}\big[1+\ln^2(y+b)+\ln^2(|x-y|)\big]~dy~\leq~\O(1),
%\]
To estimate derivatives of ${\bf G}\big[\chi_{[0,\infty[}\phi_b\big](x)$, we consider two cases:

$\bullet$ If $x+b>0$ then 
\begin{multline*}
\left|{d\over dx}{\bf G}\big[\chi_{[0,\infty[}\phi_b\big](x)\right|~=~\left|\int_{|x|}^{2+|x|}\Phi'(x+b+z)\cdot K(-z)~dz\right|\\
~\leq~\O(1)\cdot \int_{|x|}^{2+|x|} {1+\big|\ln|x+b+z|\Big|\over z}~dz
~\leq~\O(1)\cdot \big(1+\big|\ln|x|\big|\big)\cdot \int_{|x|}^{2+|x|} {1\over z}~dz~\leq~\O(1)\cdot \ln^2|x|,
\end{multline*}
and
\begin{multline*}
\left|{d^2\over dx^2}{\bf G}\big[\chi_{[0,\infty[}\phi_b\big](x)\right|~=~\left|-\Phi'(x+b+z)\cdot K(-z)\Big|^{2+|x|}_{|x|}+\int_{|x|}^{2+|x|}\Phi''(x+b+z)\cdot K(-z)~dz\right|\\
~\leq~\O(1)\cdot {1\over |x|}\cdot\left(|\ln b|+\int_{|x|}^{2+|x|} {1\over z}~dz\right)~\leq~\ds\O(1)\cdot \left|{\ln |x|\over x}\right|.
\end{multline*}

$\bullet$ Otherwise, if  $x+b\leq 0$ then 
\begin{multline*}
\left|{d\over dx}{\bf G}\big[\chi_{[0,\infty[}\phi_b\big](x)\right|~=~\left|\int_{|x|}^{2+|x|}\Phi'(x+b+z)\cdot K(-z)~dz\right|\\
~\leq~\O(1)\cdot\left[\big|\ln |x|\big|+\int_{|x|}^{2+|x|}\left|\Phi(x+b+z)\cdot K'(-z)\right|~dz\right]~\leq~\O(1)\cdot \ln^2|x|,
%\\~\leq~\O(1)\cdot\left[1+\int_{|x|}^{2+|x|} {x+b+z\over z^2}\cdot |\ln (x+b+z)|~dz\right]~\leq~\O(1)\cdot \int_{|x|}^{2+|x|} {|\ln z|\over z}~dz~\leq~\O(1)\cdot \ln^2|x|,
\end{multline*}
and
\begin{eqnarray*}
\left|{d^2\over dx^2}{\bf G}\big[\chi_{[0,\infty[}\phi_b\big](x)\right|&\leq &\O(1)\cdot \left|{\ln|x|\over x}\right|+\left|\int_{|x|}^{2+|x|}\Phi'(x+b+z)\cdot K'(-z)~dz\right|\\%~\leq~\O(1)\cdot\int_{|x|}^{2+|x|} {x+z\over z^3}\cdot |\ln (x+z)|~dz\\
&\leq&\O(1)\cdot \left[ \left|{\ln |x|\over x}\right|+\int_{|x|}^{2+|x|}\left|\Phi(x+b+z)\cdot K''(-z)\right|~dz\right]~\leq~\ds\O(1)\cdot \left|{\ln |x|\over x}\right|.
\end{eqnarray*}

\noindent {\bf Case 2.} Assume that $0<x<1/4$. We write 
\begin{eqnarray*}
{\bf G}[\chi_{[0,\infty[} \phi_b] (x)&=& -\int_0^{1-b} \Lambda(y+b) \Lambda(x-y) \, dy - \int_{1-b}^2 \Phi'(y+b) \Lambda(x-y)dy~\doteq~ -I_1-I_2.
\end{eqnarray*}
Since $\Lambda$ is $C^{3}$ in  $[1/2,2]$, it holds 
\[
\left|I^{(i)}_2(x)\right|~\leq~\O(1)\qquad\forall i\in\{0,1,2\}.
\]
On the other hand, we split $I_1$ into  three parts as follows
\[
\begin{split}
I_1&=~\int_{0}^{x/2}\Lambda(y+b)\cdot\Lambda(x-y)dy+\int_{x/2}^{3x/2}\Lambda(y+b)\cdot\Lambda(x-y)dy
\\
&\qquad\qquad+\int_{3x/2}^{1-b}\Lambda(y+b)\cdot\Lambda(x-y)dy~=~I_{11}+I_{12}+I_{13}.
\end{split}
\]
%\bel{I1}
%I_1~=~\int_{0}^{x/2}\Lambda(y)\cdot\Lambda(x-y)dy+\int_{x/2}^{3x/2}\Lambda(y)\cdot\Lambda(x-y)dy
%+\int_{3x/2}^{2}\Lambda(y)\cdot\Lambda(x-y)dy~:=~I_{11}+I_{12}+I_{13}.
%\eeq
We  estimate
\[
\begin{cases}
\big|I_{11}(x)\big|~\leq~\ds \O(1)\cdot\int_{0}^{{x\over 2}}|\ln (y+b)| |\ln(x-y)| \, dy~\leq~\O(1)\cdot \big|\ln x\big| \int_{0}^{{x\over 2}}\big|\ln y\big|dy~\leq~\O(1)\cdot x\ln^2x,\\[4mm]
\big|I'_{11}(x)\big|~=~\left|\ds\int_{0}^{{x\over 2}}\Lambda(y+b) K(x-y) \, dy+{\Lambda({x\over 2}+b)\Lambda({x\over 2})\over 2} \right|~\leq~\O(1)\cdot \ln^2x,\\[4mm]
\big|I''_{11}(x)\big|~=~\left|\ds\int_{0}^{{x\over 2}}\Lambda(y+b) K'(x-y)dy+{3\Lambda({x\over 2}+b) K({x\over 2})\over 4}+{\Lambda({x\over 2})K({x\over 2}+b)\over 4} \right|~\leq~\O(1)\cdot \left|\ds{\ln x\over x}\right|,
\end{cases}
\]
and 
\[
\begin{cases}
\big|I_{13}(x)\big|&\leq~\ds \O(1)\cdot\int_{{3x\over 2}}^{1-b}(1+|\ln (y+b)|)\cdot (1+|\ln(y-x)|) \, dy~\leq~\O(1) \cdot x\ln^2x,\\[4mm]
\big|I'_{13}(x)\big|&=~\left|\ds\int_{{3x\over 2}}^{1-b}\Lambda(y+b) K(x-y) \, dy - {3\Lambda({3x\over 2}+b) \Lambda(-{x\over 2})\over 2} \right|~\leq~\O(1)\cdot \ln^2x,\\[5mm]
\big|I''_{13}(x)\big|&=~\left|\ds\int_{{3x\over 2}}^{1-b}\Lambda(y+b) K'(x-y) \, dy-{9K({3x\over 2}+b) \Lambda(-{x\over 2})\over 4} - {3\Lambda({3x\over 2}+b)K(-{x\over 2})\over 4} \right|\\&\leq~\O(1)\cdot \left|\ds{\ln x\over x}\right|.
\end{cases}
\]
Concerning $I_{12}$, first of all, by a change of variable, it holds
\begin{eqnarray*}
I_{12}(x) = \lim_{\ve\to 0+} \left(\int_{x/2}^{x-\ve}+\int_{x+\ve}^{{3x\over 2}}\right)\Lambda(y+b)\cdot \Lambda(x-y) \, dy = \int_{-{x\over 2}}^{{x\over 2}} \Lambda(x+b-z) \cdot \Lambda(z) \, dz,
\end{eqnarray*}
and one directly computes that 
\[
\begin{cases}
\big|I_{12}(x)\big| & \leq ~ \ds \O(1) \cdot \big|\ln x\big| \cdot \int_{0}^{{x\over 2}}|\ln z| \, dz \leq \O(1)\cdot x\ln^2x,\\[4mm]
\big|I'_{12}(x)\big| & = ~  \left|{\ds \int_{-{x\over 2}}^{{x\over 2}}} K(x+b-z) \cdot \Lambda(z) \, dz + {1\over 2}\cdot \left[\Lambda({3x\over 2}+b)\Lambda(-{x\over 2})+\Lambda({x\over 2}+b)\Lambda({x\over 2})\right]\right|\\[4mm]
& \leq ~ \ds \O(1) \cdot \left({1\over x}\cdot \int_{0}^{{x\over 2}}|\ln z| \, dz + \ln^2x\right) \leq \O(1) \cdot \ln^2x,\\[4mm]
\big|I''_{12}(x)\big| & \leq ~ \O(1) \cdot \left(\left| \ds \int_{-{x\over 2}}^{{x\over 2}} K'(x+b-z) \cdot \Lambda(z) \, dz\right| + \left|\ds{\ln x\over x}\right|\right)\\[4mm]
& \leq ~ \O(1) \cdot \left(\ds {1\over x^2} \cdot \int_{0}^{{x\over 2}}\big|\ln z\big| \, dz +\left|\ds{\ln x\over x}\right|\right) \leq \O(1)\cdot \left|\ds{\ln x\over x}\right|.
\end{cases}
\]
From the previous step, we obtain the first three estimates in (\ref{Gp-1}) for  $0<|x|<1/4$.
\medskip

Finally, observe that $\chi_{[0,\infty[}\phi_b$ is continuous with compact support and smooth outside the origin. Hence, ${\bf G}\big[\chi_{[0,\infty[}\phi_b\big]$ is smooth outside the origin. As $|x|\to+\infty$, for $i\in\{1,2\}$ we have 
\[
\left|{\bf G}\big[\chi_{[0,\infty[}\phi_b\big](x)\right|~\leq~\O(1)\cdot x^{-1},\qquad \left|{d^i\over dx^i}{\bf G}\big[\chi_{[0,\infty[}\phi_b\big](x)\right|~\leq~\O(1)\cdot x^{-(1+i)},
\]
%As a matter of fact, for every $|x| > 3$ and every sufficiently small $\ve > 0$ there holds
%\[
%\left\{y \in \R : |y-x| > \ve\right\} \cap \left\{y \in \R : |y| < 2\right\} = \left\{y \in \R : |y| < 2\right\},
%\]
%hence
%\[
%{\bf G}\big[\chi_{[0,\infty[}\Phi\big](x) = \lim_{\ve \to 0+} \int_{\{|y-x| > \ve\} \cap \{|y| < 2\}} K(x-y) \chi_{[0,\infty[}(y) \Phi(y) \, dy = \int_0^2 K(x-y) \Phi(y) \, dy
%\]
%and one obtains the estimates above.
and using the first three estimates in (\ref{Gp-1}) we obtain the last estimate in (\ref{Gp-1}).
\endproof

Following the same argument in Lemma \ref{A1}, one can show that

\begin{remark}\label{Ln}   Given $\lambda_1,\lambda_2\in \R$, the  function $$v(x)= \left(\lambda_1\cdot\chi_{]-\infty,0[}+ \lambda_2\cdot\chi_{]0,\infty[}\right)\cdot \eta(x)x$$  is more regular than $\Phi$. Thus, one can follow the same argument as in Lemma \ref{A1} to obtain  for all $|x|<1/2$ that 
\[
\left|{\bf G}[v](x)\right|~\leq~\O(1)\cdot \big(|\lambda_1|+|\lambda_2|\big),\qquad \left|{d\over dx}{\bf G}[v](x)\right|~\leq~\O(1)\cdot \big(|\lambda_1|+|\lambda_2|\big)  \cdot \ln^2|x|,
\]
and 
\[
\big\|{\bf G}[v]\big\|_{H^2(\R\backslash [-\delta,\delta])}~\leq~\O(1)\cdot \big(|\lambda_1|+|\lambda_2|\big)\cdot \delta^{-2/3}\qquad\forall \delta>0.
\]
\end{remark}

\begin{lemma}\label{B-est} Let $v\in H^{2}(\R\backslash\{0\})$ be such that 
\[
\|v\|_{H^{i}(\R\backslash\{0\})}~\leq~M_i,\qquad i\in \{1,2\}.
\]
Set $D^{(v)}(x)\doteq{\bf G}[v](x)-\big[{v(0+)-v(0-)}\big]\cdot \eta(x)\cdot \Lambda(x)$. Then for every $0<|x|<1/2$ and $\delta>0$ small, we have that  
\bel{B1}
\begin{cases}
\left|D^{(v)}(x)\right|~\leq~\O(1)\cdot M_1,\qquad \left|D^{(v)}_x(x)\right|~\leq~\O(1)\cdot M_1\cdot \ln^2|x|,\\[4mm]
\ds \big\|D^{(v)}\big\|_{H^1(\R\backslash \{0\})}~\leq~\O(1)\cdot {M_1},\qquad \big\|D^{(v)}\big\|_{H^2(\R\backslash [-\delta,\delta])}~\leq~\O(1)\cdot M_2\cdot\delta^{-2/3}.
\end{cases}
\eeq
Consequently, ${\bf G}[v]$ is in $H^1_{loc}(\R\backslash \{0\})$ and
\[
\|{\bf G}[v]\|_{H^1(\R\backslash [-\delta,\delta])}~\leq~M_1\cdot \delta^{-1/2}.
\]
\end{lemma}
{\bf Proof.}   We first split $v$ into two parts 
\[
v=v_1+v_2,\qquad \text{where} \qquad v_1(x)~=~\begin{cases}
v(0-)\cdot \eta(x)&\mathrm{if}~~ x<0,\\[3mm]
v(0+)\cdot \eta(x)&\mathrm{if}~~ x>0.
\end{cases}
\]
For $i\in \{1,2\}$, one has that 
\[
\|v_1\|_{H^{(i)}(\R\backslash\{0\})}~\leq~\big(|v(0-)|+|v(0+)|\big)\cdot \|\eta\|_{H^{(i)}(\R\backslash\{0\})}~\leq~4M_1\cdot \|\eta\|_{H^{(i)}(\R\backslash\{0\})},
\]
and the function $v_2\in H^1(\R)$ satisfies 
\[
\|v_2\|_{H^{1}(\R)}~\leq~\left(1+4\|\eta\|_{H^{1}(\R)}\right)\cdot M_1,\qquad \|v_2\|_{H^{2}(\R\backslash\{0\})}~\leq~\left(1+4\|\eta\|_{H^{2}(\R)}\right)\cdot M_2.
\]

{\bf Step 1.} For every $0<|x|<1/2$, we have
\bel{Gv1}
{\bf G}[v_1](x)=\int_{-2}^{2}v_1'(y)\cdot \Lambda(x-y)~dy+\big[v(0+)-v(0-)\big]\cdot \Lambda(x)=I_1+\big[v(0+)-v(0-)\big]\cdot \Lambda(x).
\eeq
Recalling (\ref{eta}), we estimate
\begin{eqnarray*}
\big|I_1(x)\big|&=&\left|v(0-)\cdot \int_{-2}^{-1}\eta'(y)\cdot \Lambda(x-y)~dy+v(0+)\cdot \int_{1}^{2}\eta'(y)\cdot \Lambda(x-y)~dy\right|\\
&\leq& \O(1)\cdot M_1\cdot \left(\int_{-2}^{-1}| \Lambda(x-y)|~dy+\int_{1}^{2}| \Lambda(x-y)|~dy\right)~\leq~ \O(1)\cdot M_1,
\end{eqnarray*}
and 
\begin{eqnarray*}
\left|I^{(i)}_1(x)\right|~\leq~ \O(1)\cdot M_1\cdot \left(\int_{-2}^{-1}\left|\Lambda^{(i)}(x-y)\right|~dy+\int_{1}^{2}\left| \Lambda^{(i)}(x-y)\right|~dy\right)~\leq~ \O(1)\cdot M_1
\end{eqnarray*}
for $i\in \{1,2\}$. Hence, (\ref{Gv1}) and (\ref{K1}) yield
\bel{G-v1-e1}
\begin{cases}
\ds {\bf G}[v_1](x)~=~\O(1)\cdot M_1+\big[v(0+)-v(0-)\big]\cdot \Lambda(x)\\[4mm]
{\bf G}'[v_1](x)~=~\O(1)\cdot M_1+\big[v(0+)-v(0-)\big]\cdot K(x),\\[4mm]
 {\bf G}''[v_1](x)~=~\O(1)\cdot M_1+\big[v(0+)-v(0-)\big]\cdot K'(x).
\end{cases}
\eeq
Moreover, since $v_1$ is bounded and compactly supported, as $|x| \to +\infty$ 
\begin{equation*}
\left|{\bf G}[v_1](x)\right|~\leq~\O(1)\cdot M_1\cdot  x^{-1},\qquad \left|{d^i\over dx^i}{\bf G}[v_1](x)\right|~\leq~\O(1)\cdot M_1\cdot x^{-(1+i)}\quad i\in\{1,2\},
\end{equation*}
and we have
\bel{es-v1}
\begin{cases}
\|{\bf G}[v_1]-\big[v(0+)-v(0-)\big]\cdot \Lambda \cdot \eta\|_{H^1(\R\backslash\{0\})}~\leq~\O(1)\cdot M_1,\\[3mm]
\|{\bf G}[v_1]-\big[v(0+)-v(0-)\big]\cdot \Lambda \cdot \eta \|_{H^2(\R\backslash\{0\})}~\leq~\O(1)\cdot M_1.
\end{cases}
\eeq

{\bf Step 2.} To estimate ${\bf G}[v_2]$,  we first recall that $v_2\in H^1(\R)$. By the continuity of the linear operator ${\bf G}:{\bf L}^2(\R)\to {\bf L}^2(\R)$, we have 
\[
\|{\bf G}[v_2]\|_{H^1(\R)}~\leq~\O(1)\cdot \|v_2\|_{H^1(\R)}~\leq~\O(1)\cdot M_1,
\]
and this particularly  yields 
\[
\big|{\bf G}[v_2](x)\big|~\leq~2\cdot \big\|{\bf G}[v_2]\big\|_{H^1(\R)}~\leq~\O(1)\cdot M_1\qquad\forall x\in \R.
\]
Thus, recalling (\ref{G-v1-e1}) and (\ref{es-v1}), we get the first and the third estimate in (\ref{B1}).
\medskip

{\bf Step 3.} To achieve the second and the fourth estimate in (\ref{B1}), we split $v_2$ into two parts 
\[
v_{2}=v_{21}+v_{22},\qquad\qquad v_{21}(x)~=~\begin{cases}
v_{2,x}(0-)\cdot x\eta(x)&\mathrm{if}~~ x<0,\\[3mm]
v_{2,x}(0+)\cdot x\eta(x)&\mathrm{if}~~ x>0.
\end{cases}
\]
Since $|v_{2,x}(0\pm)|\leq 2\cdot \|v_{2,x}\|_{H^1(\R\backslash\{0\})}\leq 2M_2$, one has that  $\|v_{22}(\cdot)\|_{H^2(\R)}\leq\O(1)\cdot M_2$. Thus, by the continuity of the linear operator ${\bf G}:{\bf L}^2(\R)\to {\bf L}^2(\R)$, we get
\[
\big\|{\bf G}[v_{22}](\cdot)\big\|_{H^2(\R)}~\leq~\O(1)\cdot M_2,\ |{\bf G}[v_{22}](x)|,\left|{d\over dx}{\bf G}[v_{22}](x)\right|~\leq~\O(1)\cdot M_2, \ \forall x\in \R.
\]
Finally, by Remark \ref{Ln} we have 
\[
\left|{\bf G}[v_{21}](x)\right|~\leq~\O(1)\cdot M_2,\qquad \left|{d\over dx}{\bf G}[v_{21}](x)\right|~\leq~\O(1)\cdot M_2  \cdot \ln^2|x|\qquad \forall |x|<1/2,
\]
and 
\[
\big\|{\bf G}[v_{21}]\big\|_{H^2(\R\backslash [-\delta,\delta])}~\leq~\O(1)\cdot M_2\cdot \delta^{-2/3}\qquad\forall \delta>0.
\]
Thus,  (\ref{G-v1-e1}) and (\ref{es-v1}) yield  the second and the fourth estimates  in (\ref{B1}).
\endproof

\v

{\bf Acknowledgments.}  Khai T. Nguyen was partially supported by National Science Foundation grant DMS-2154201. 
J.Schino is a member of GNAMPA (INdAM). Part of this work was done when J. Schino was an employee of North Carolina State University, to which he expresses deep gratitude.

\end{document}